\documentclass[english]{article}
\usepackage[T1]{fontenc}
\usepackage[latin9]{inputenc}
\usepackage[a4paper]{geometry}
\usepackage{color}
\usepackage{babel}
\usepackage{longtable}
\usepackage{float}
\usepackage{booktabs}
\usepackage{textcomp}
\usepackage{enumitem}
\usepackage{amsmath}
\usepackage{amsthm}
\usepackage{amssymb}
\usepackage{graphicx}
\usepackage[authoryear]{natbib}
\usepackage[unicode=true,pdfusetitle,
 bookmarks=true,bookmarksnumbered=false,bookmarksopen=false,
 breaklinks=false,pdfborder={0 0 1},backref=false,colorlinks=false]
 {hyperref}

\makeatletter

\providecommand{\tabularnewline}{\\}
\floatstyle{ruled}
\newfloat{algorithm}{tbp}{loa}
\providecommand{\algorithmname}{Algorithm}
\floatname{algorithm}{\protect\algorithmname}

\numberwithin{equation}{section}
\theoremstyle{plain}
\newtheorem{thm}{\protect\theoremname}[section]
\theoremstyle{plain}
\newtheorem{prop}[thm]{\protect\propositionname}

    \usepackage{graphicx}

    \usepackage{enumitem}

    \setlist{itemsep=0pt, parsep = 0pt}

    \setlist[enumerate,1]{label=(\roman*), ref=(\roman*) }
    \setlist[enumerate,2]{label=(\alph*),  ref=(\alph*) }
 
\usepackage{setspace}
\AtBeginDocument{%
	\addtolength\abovedisplayskip{-0.4\baselineskip}%
	\addtolength\belowdisplayskip{-0.4\baselineskip}%
}

\AtBeginDocument{

}

\makeatother

\providecommand{\propositionname}{Proposition}
\providecommand{\theoremname}{Theorem}

\title{Numerical valuation of European options under two-asset infinite-activity
exponential L{\'e}vy models}
\author{
Massimiliano Moda\thanks{Department of Mathematics, University of Antwerp, Middelheimlaan 1, 2020 Antwerp, Belgium.\hfill \break Email: massimiliano.moda@uantwerpen.be},\, 
Karel J. in 't Hout\thanks{Department of Mathematics, University of Antwerp, Middelheimlaan 1, 2020 Antwerp, Belgium.\hfill\break Email: karel.inthout@uantwerpen.be},\,
Mich{\`e}le Vanmaele\thanks{Department of Mathematics, Computer Science and Statistics, Ghent University, 9000 Ghent, Belgium.\hfill\break Email: michele.vanmaele@ugent.be},\,
Fred Espen Benth\thanks{Department of Data Science and Analytics, BI Norwegian Business School, N-0484 Oslo, Norway.\hfill\break Email: fred.e.benth@bi.no}
}

\begin{document}

\maketitle


\begin{abstract}
We propose a numerical method for the valuation of European-style options under two-asset infinite-activity exponential L{\'e}vy models. Our method extends the effective approach developed by \citet{wang2007} for the 1-dimensional case to the 2-dimensional setting and is applicable for general L{\'e}vy measures under mild assumptions. A tailored discretization of the non-local integral term is developed, which can be efficiently evaluated by means of the fast Fourier transform. For the temporal discretization, the semi-Lagrangian $\theta$-method is employed in a convenient splitting fashion, where the diffusion term is treated implicitly and the integral term is handled explicitly by a fixed-point iteration. Numerical experiments for put-on-the-average options under Normal Tempered Stable dynamics reveal a favourable convergence behaviour of our method whenever the exponential L{\'e}vy process has finite-variation. 
In addition, a relevant theoretical convergence result for the discretization of the integral term is proved.
\end{abstract}

\section{Introduction}

The accurate valuation of derivative securities in modern financial
markets requires modeling techniques capable of capturing empirical
irregularities in asset price dynamics. Classical models based on
Brownian motion, such as the Black--Scholes model, rely on continuous-path
diffusion and fail to reflect important stylized facts, such as heavy
tails and skewness in log-returns. This has motivated the use of L{\'e}vy
processes in the last decades, which naturally offer a richer class
of models for asset dynamics. Among various L{\'e}vy models, the Normal
Inverse Gaussian (NIG) process has emerged as a parsimonious and effective
choice to capture such characteristics. Among others, \citet{rydberg1997}
shows how the NIG model provides a significantly better statistical
fit to equity return data compared to classical Gaussian-based models. L{\'e}vy
models allow for a more realistic representation of market risk and
are therefore natural candidates for use in option pricing models.

In this paper, we propose a numerical method for pricing European-style financial
derivatives written on two underlying assets, whose dynamics are driven
by a 2-dimensional L{\'e}vy process, with particular focus on infinite-activity processes. 
Financial pricing under jump-diffusion models
can be approached through various methodologies, such as Monte Carlo
simulation, Fourier-based methods, and partial integro-differential
equations (PIDEs). Monte Carlo methods are flexible and easy to implement,
but they suffer from slow convergence. Fourier-based methods, such as
in \citet{jackson2008} and \citet{ruijter2012}, can be applied when
the characteristic exponent of the process is known, and they can achieve
exponential convergence. Numerical methods for PIDEs, such as in \citet{cont2005},
\citet{dhalluin2005}, \citet{wang2007}, \cite{clift2008}, \cite{salmi2014}, 
\cite{kaushansky2018}, \citet{boen2021} and \citet{inthout2023},
can instead be applied when the L{\'e}vy measure is known, and do not
require knowledge of the characteristic exponent.
They are applicable to a wide variety of financial derivatives.

The numerical method derived in this paper focuses on the case where the underlying
2-dimensional L{\'e}vy process exhibits infinite-activity, meaning that
an infinite number of jumps occur over any finite time horizon. In
this setting, particular care must be taken in the discretization
of the non-local 2-dimensional integral term in the PIDE near
the origin, where the L{\'e}vy measure becomes singular. 

The main contribution of this paper is an extension of the effective numerical 
solution approach of \citet{wang2007} from the 1-dimensional to the 2-dimensional 
setting. 
Here, a key idea, originally introduced in \citet{asmussen2001} and \citet{cont2005}, 
is to replace the small jumps with an artificial diffusion term. This substitution 
enables the development of a tailored quadrature scheme.
For the efficient evaluation of the discretized integral operator, a fast Fourier 
transform (FFT) algorithm is constructed.
For the temporal discretization, the semi-Lagrangian $\theta$-method is considered.
Here, operator splitting is applied, where the diffusion term is treated implicitly
and the integral term is handled explicitly by a fixed-point iteration.
For the large linear system in each time step, the BiCGSTAB iterative solver is used.

To assess the performance of the proposed numerical method, we conduct numerical experiments for put-on-the-average options under Normal Tempered Stable dynamics. These experiments reveal a favourable convergence behaviour whenever the underlying exponential L{\'e}vy process has finite-variation. A relevant theoretical result is proved on the order of convergence for the discretization of the integral term.

An outline of this paper is as follows. In Section~\ref{sec: model framework},
we introduce the market model and the PIDE for the derivative pricing.
In Section~\ref{sec: Numerical scheme} the proposed numerical scheme is derived
and a theorem on the approximation error for the integral term is presented.
Numerical experiments are discussed in Section~\ref{sec: Numerical experiment}.
The final Section~\ref{sec: Conclusions} contains our conclusions.

\section{\label{sec: model framework}Model framework}

\subsection{Market model}

Let $\left(\Omega,\mathcal{F},\left(\mathcal{F}_{t}\right)_{t\in\left[0,T\right]},\mathbb{P}\right)$
be a filtered probability space, for some given $T>0$. We consider
an arbitrage free market characterized by a constant (instantaneous)
risk-free interest rate $r$ and an equivalent martingale measure
$\mathbb{Q}\sim\mathbb{P}$. We assume that there exist two risky
assets whose prices are modeled by the 2-dimensional stochastic process
$X=\left(X^{\left(1\right)},X^{\left(2\right)}\right)$ that solves
the following stochastic differential equation
\begin{equation}
dX\left(t\right)=\mu\left(X\left(t\right)\right)dt+\Sigma\left(X\left(t\right)\right)dW\left(t\right)+\int_{\mathbb{R}_{*}^{2}}\gamma\left(z,X\left(t_{-}\right)\right)\tilde{\Pi}\left(dt,dz\right)\qquad\left(t\in\left(0,T\right]\right)\label{eq: intro - asset prices dynamics}
\end{equation}
for some non-negative initial value $X\left(0\right)$. In (\ref{eq: intro - asset prices dynamics}),
$W$ denotes a standard 2-dimensional Wiener process and $\tilde{\Pi}$ is
a compensated Poisson random measure with L{\'e}vy measure $\ell$ over
$\mathbb{R}_{*}^{2}=\mathbb{R}^{2}\setminus\left\{ 0\right\} $. Both
are directly defined under $\mathbb{Q}$ and are mutually independent. 

The functions $\mu:\mathbb{R}_{\geq0}^{2}\rightarrow\mathbb{R}^{2}$,
$\Sigma:\mathbb{R}_{\geq0}^{2}\rightarrow\mathbb{R}^{2\times2}$ and
$\gamma:\mathbb{R}^{2}\times\mathbb{R}_{\geq0}^{2}\rightarrow\mathbb{R}^{2}$
are called drift, diffusion, and jump function (or term) respectively,
where $\mathbb{R}_{\geq0}^{2}=\left\{ x\in\mathbb{R}^{2}:x^{\left(i\right)}\geq0\text{ for }i=1,2\right\} $.
In this paper, we consider the case of the well-known exponential
L{\'e}vy process, i.e. where the coordinates of the coefficient functions
are defined for $i,j=1,2$ as follows 
\begin{align}
\mu^{\left(i\right)}\left(x\right) & =x^{\left(i\right)}r\label{eq: intro - drift function}\\
\left(\Sigma\Sigma^{\top}\right)^{\left(i,j\right)}\left(x\right) & =x^{\left(i\right)}x^{\left(j\right)}\left(\sigma\sigma^{\top}\right)^{\left(i,j\right)}\label{eq: intro - diffusion function}\\
\gamma^{\left(i\right)}\left(z,x\right) & =x^{\left(i\right)}\left(e^{z^{\left(i\right)}}-1\right),\label{eq: intro - jump function}
\end{align}
where $\sigma\sigma^{\top}$ is a constant positive
definite symmetric $2\times2$ matrix and $\sigma$ denotes the volatility matrix.
Here, $\Sigma\Sigma^{\top}\left(x\right)$ is a shorthand notation
for the matrix product $\Sigma\left(x\right)\Sigma^{\top}\left(x\right)$. 

Let $\left\Vert \cdot\right\Vert $ be any given norm on $\mathbb{R}^{2}$.
In this paper we assume that $\ell$ is absolutely continuous, the corresponding L{\'e}vy process has finite variance, i.e. 
\begin{equation}\label{eq: intro - finite variance hypothesis}
\int_{\mathbb{R}^2_*}\left\Vert z\right\Vert^2\ell\left(dz\right)<\infty,
\end{equation}
 and there exist constants $A_{\ell}<2$, $B_{\ell}>2$ and $C_{\ell}\left(h\right)$, $C_{\ell}^{\prime}\left(h\right)>0$, for any $h>0$, such that
\begin{equation}\label{eq: intro - upper bound of the L=0000E9vy measure}
\begin{cases}
\ell\left(z\right)\leq C_{\ell}\left(h\right)\left\Vert z\right\Vert^{-2-A_{\ell}} & \text{for any }z\text{ such that }\left\Vert z\right\Vert\in\left(0,h\right],\\

\left|\ell_{z}^{\left(j\right)}\left(z\right)\right|\leq C_{\ell}^{\prime}\left(h\right)\left\Vert z\right\Vert ^{-3-A_{\ell}}& \text{for any }z\text{ such that }\left\Vert z\right\Vert\in\left(0,h\right],\\

\ell\left(z\right)=O\left(e^{-B_{\ell}\left\Vert z\right\Vert }\right) & \text{as }\left\Vert z\right\Vert \rightarrow\infty,
\end{cases}
\end{equation} where $\ell_{z}^{\left(j\right)}$ denotes the $j$-th partial derivative of $\ell$ with respect to $z$, with $j=1,2$.
The number $A_{\ell}$ governs the activity and variation of the
process: $X$ is of finite-activity if $A_{\ell}<0$, since $\int_{\mathbb{R}^{2}_*}\ell\left(dz\right)<\infty$;
it is of finite-variation if $A_{\ell}<1$, since $\int_{0 < \left\Vert z\right\Vert <\epsilon}\left\Vert z\right\Vert \ell\left(dz\right)<\infty$
for any $\epsilon>0$. The number $B_\ell$ characterizes the exponential decay of $\ell$  at infinity. Since the process $X$ has finite moments of all orders up to $k\in\mathbb{N}$  if and only if  $\int_{\left\Vert z\right\Vert >\epsilon}e^{k\left\Vert z\right\Vert }\ell\left(dz\right)<\infty$ for any $\epsilon>0$, then $k<B_\ell$  provides a necessary condition of it. Following \citet[Chapter 6]{applebaum2004}, the stronger condition $B_{\ell}\geq2$ is necessary to guarantee the
existence of a unique solution with finite variance to the stochastic
differential equation (\ref{eq: intro - asset prices dynamics}).
Most of the common L{\'e}vy processes in finance satisfy the conditions 
(\ref{eq: intro - upper bound of the L=0000E9vy measure}), such as
the Kou, Carr--Geman--Madan--Yor (CGMY), Variance Gamma (VG) and Normal
Inverse Gaussian (NIG) models.

In this work, we focus on the case of 2-dimensional Normal Tempered
Stable (NTS) processes. These are obtained by subordinating a 2-dimensional
arithmetic Brownian motion with a Tempered Stable subordinator. A
detailed construction of the NTS process together with its main properties
is provided in Appendix \ref{app: Normal Tempered Stable process}.
The choice of this class of processes is motivated by two reasons.
First, bivariate VG and NIG processes arise as particular cases. Second,
the associated L{\'e}vy measure satisfies the conditions (\ref{eq: intro - upper bound of the L=0000E9vy measure})
with constant $A_{\ell}=2\alpha$, where $\alpha$ is the key model
parameter. The NTS framework provides a convenient and flexible setting
for the purposes of this paper. 

\subsection{Initial boundary value problem for derivatives pricing}

By the fundamental theorem of asset-pricing, the value at time $t\in\left[0,T\right]$
of an European-style\footnote{Means a financial derivative with no intermediate cash flows.}
financial derivative of $X$ with maturity $T$ is represented by
the stochastic process $P$ such that
\[
P\left(t\right)=\mathbb{E}^{\mathbb{Q}}\left[\phi\left(X\left(T\right)\right)e^{-r\left(T-t\right)}\mid\mathcal{F}_{t}\right]
\]
where $\phi:\mathbb{R}^{2}\rightarrow\mathbb{R}$ denotes the pay-off
function and $\mathbb{E}^{\mathbb{Q}}\left[\cdot\mid\mathcal{F}_{t}\right]$
is the $\mathcal{F}_{t}$-conditional expected value (i.e. knowing
the history of the asset prices up to $t$) under $\mathbb{Q}$.

Let $\mathcal{A}$ be the infinitesimal generator of $X$ (see \citet{applebaum2004},
\citet{garroni1992} and \citet{oksendal2019}), defined in matrix
notation as\footnote{By expanding the term, we obtain the common notation used for $\mathcal{A}$,
which is 
\begin{align*}
\mathcal{A}u\left(x,t\right)= & \sum_{i=1}^{2}\mu^{\left(i\right)}\left(x\right)\frac{\partial u}{\partial x^{\left(i\right)}}\left(x,t\right)+\frac{1}{2}\sum_{i,j=1}^{2}\left(\Sigma\Sigma^{\top}\right)^{\left(i,j\right)}\left(x\right)\frac{\partial^{2}u}{\partial x^{\left(i\right)}\partial x^{\left(j\right)}}\left(x,t\right)\\
 & +\int_{\mathbb{R}_{*}^{2}}\left(u\left(x+\gamma\left(z,x\right),t\right)-u\left(x,t\right)-\sum_{i=1}^{2}\gamma^{\left(i\right)}\left(z,x\right)\frac{\partial u}{\partial x^{\left(i\right)}}\left(x,t\right)\right)\ell\left(dz\right).
\end{align*}
}
\begin{equation}
\mathcal{A}u\left(x,t\right)=\mu\left(x\right)^{\top}u_{x}\left(x,t\right)+\frac{1}{2}\mathbf{1}^{\top}\left(u_{xx}\left(x,t\right)\circ\Sigma\Sigma^{\top}\left(x\right)\right)\mathbf{1}+\int_{\mathbb{R}_{*}^{2}}f\left(z,x,t\right)\ell\left(dz\right)\label{eq: intro - infinitesimal generator, definition}
\end{equation}
where $\mathbf{1}=\left[1,1\right]^{\top}$, the symbol $\circ$ denotes
the Hadamard (element-wise) product\footnote{In this paper, we use the convention $AB\circ CD=\left(AB\right)\circ\left(CD\right)$,
for any suitable matrices $A,B,C,D$.} and 
\begin{equation}
f\left(z,x,t\right)=u\left(x+\gamma\left(z,x\right),t\right)-u\left(x,t\right)-\gamma\left(z,x\right)^{\top}u_{x}\left(x,t\right).\label{eq: intro - integrand function}
\end{equation}
If there exists a function $u:\mathbb{R}_{\geq0}^{2}\times\left[0,T\right]\rightarrow\mathbb{R}$
that solves the following initial value problem for a partial integro-differential
equation (PIDE)
\begin{equation}
\begin{cases}
u_{t}\left(x,t\right)=\mathcal{A}u\left(x,t\right)-ru\left(x,t\right) & \text{for any }\left(x,t\right)\in\mathbb{R}_{\geq0}^{2}\times\left(0,T\right]\\
u\left(x,0\right)=\phi\left(x\right)
\end{cases}\label{eq: intro - IVP (1)}
\end{equation}
then the value of the financial derivative is given by
\[
P\left(t\right)=u\left(X\left(t\right),T-t\right).
\]

Note that $u$ also satisfies the PIDE at the boundary of the $x$-domain,
as in the case of option pricing with the Black--Scholes model.

\section{\label{sec: Numerical scheme}Numerical scheme}

In this section, we derive the numerical scheme proposed for problem
(\ref{eq: intro - IVP (1)}).

The method consists of three main steps: integral discretization,
spatial discretization, and temporal discretization. By discretization,
we mean that the pertinent integro/differential operators are approximated
on a given finite set of grid points. The adjectives indicate the
variable being discretized: integral for $z$, spatial for $x$, and
temporal for $t$.

The integral discretization yields an approximation to the integral
term in (\ref{eq: intro - IVP (1)}) for any given pair $\left(x,t\right)\in\mathbb{R}_{\geq0}^{2}\times\left[0,T\right]$.
The quadrature formula that we derive is inspired by the ideas in
\citet{wang2007} and \citet{cont2005}, where the singular part of
the integral near the origin $z=0$ is approximated by a diffusion
(second-order spatial derivative). The integral discretization leads
to the approximate problem (\ref{eq: integral disc - IVP (2)}) where
the integral in (\ref{eq: intro - IVP (1)}) has been replaced by
a summation term.

The spatial discretization concerns the diffusion and summation terms
in (\ref{eq: integral disc - IVP (2)}). For the diffusion term, a
standard second-order central finite difference scheme is applied
on a suitable nonuniform spatial grid. For the summation term, a direct
valuation on the spatial grid is computationally too expensive. For
the efficient treatment of this term, we shall extend the FFT-based
approach by \citet{wang2007}.

The temporal discretization is done by the semi-Lagrangian $\theta$-method.
The semi-Lagrangian approach is chosen to take into account that problem
(\ref{eq: integral disc - IVP (2)}) can be convection-dominated.
In each time step of the semi-Lagrangian $\theta$-method, the summation
term appears in an implicit manner. To effectively handle this term,
a fixed-point iteration is employed.

\subsection{\label{subsec: integral disc}Integral discretization}

When the L{\'e}vy measure is singular, classical quadrature formulae such as the 
midpoint rule or the trapezoidal rule fail.
In fact, in this case the error will blow up as the number of quadrature
points increases. To address this problem, we develop in this subsection
a different quadrature formula.

First, it is useful to investigate the behaviour of $f$, defined
in (\ref{eq: intro - integrand function}), around the origin with
respect to $z$. For any given point $\left(x,t\right)\in\mathbb{R}_{\geq0}^{2}\times\left[0,T\right]$,
the Taylor approximation of the function $z\mapsto f\left(z,x,t\right)$
at $z=0$ is given by 
\[
f\left(z,x,t\right)=f\left(0,x,t\right)+z^{\top}f_{z}\left(0,x,t\right)+\frac{1}{2}z^{\top}f_{zz}\left(0,x,t\right)z+\varepsilon\left(z,x,t\right)\qquad\text{as }\left\Vert z\right\Vert \rightarrow0^{+},
\]
where $f_{z}$ and $f_{zz}$ are the gradient and the Hessian of $f$
with respect to $z$. Here, $\varepsilon$ denotes the remainder and
is such that $\varepsilon\left(z,x,t\right)=O\left(\left\Vert z\right\Vert ^{3}\right)$.
Invoking (\ref{eq: intro - integrand function}) and noting that $f\left(0,x,t\right)=0$
and $f_{z}\left(0,x,t\right)=0$, we can rewrite the previous equation,
after some straightforward computations, as
\begin{equation}
f\left(z,x,t\right)=\frac{1}{2}\mathbf{1}^{\top}\left(u_{xx}\left(x,t\right)\circ I_{x}zz^{\top}I_{x}\right)\mathbf{1}+\varepsilon\left(z,x,t\right)\qquad\text{as }\left\Vert z\right\Vert \rightarrow0^{+},\label{eq: integral disc - Taylor approximation of the integrand}
\end{equation}
where $u_{xx}$ is the Hessian of $u$ with respect to $x$ and $I_{x}=\text{diag}\left(x^{\left(1\right)},x^{\left(2\right)}\right)$.

\begin{figure}
\caption{\label{fig: integral disc - Partition of the integration domain}Partition
of the integration domain $R_{z}$}

\noindent \centering{}\includegraphics{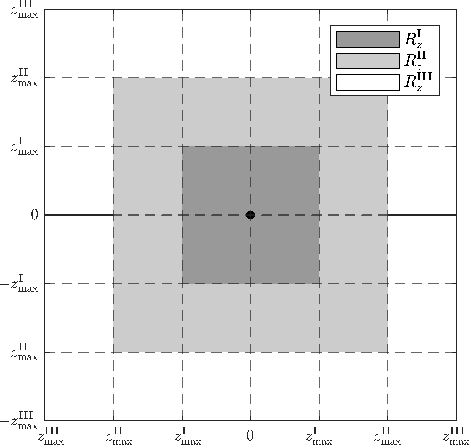}
\end{figure}

Next, let $R_{z}^{\mathbf{I}}$, $R_{z}^{\mathbf{II}}$ and $R_{z}^{\mathbf{III}}$
be three sets defined by
\begin{align*}
R_{z}^{\mathbf{I}} & =\left\{ z\in\mathbb{R}^{2}:\left\Vert z\right\Vert _{\infty}\leq z_{\max}^{\mathbf{I}}\right\} ,\\
R_{z}^{\mathbf{II}} & =\left\{ z\in\mathbb{R}^{2}:z_{\max}^{\mathbf{I}}<\left\Vert z\right\Vert _{\infty}\leq z_{\max}^{\mathbf{II}}\right\} ,\\
R_{z}^{\mathbf{III}} & =\left\{ z\in\mathbb{R}^{2}:z_{\max}^{\mathbf{II}}<\left\Vert z\right\Vert _{\infty}\leq z_{\max}^{\mathbf{III}}\right\} ,
\end{align*}
where $\left\Vert z\right\Vert _{\infty}=\max_{j=1,2}\left|z^{\left(j\right)}\right|$
and $0<z_{\max}^{\mathbf{I}}<z_{\max}^{\mathbf{II}}<z_{\max}^{\mathbf{III}}$
are given numbers, which will be specified below.
The above three sets represent a partition of $R_{z}  =\left\{ z\in\mathbb{R}^{2}: \left\Vert z\right\Vert _{\infty}\leq z_{\max}^{\mathbf{III}}\right\}$, which is a square centered
at the origin, as shown in Figure \ref{fig: integral disc - Partition of the integration domain}.
 For any given $N_{z}\in\mathbb{N}$, define a set of points $\mathbf{z}$
whose elements are 
\[
z_{l_{1}l_{2}}=\left(\left(l_{1}+\frac{1}{2}\right)h_{z},\left(l_{2}+\frac{1}{2}\right)h_{z}\right)\qquad\left(l_{1},l_{2}=-N_{z},-N_{z}+1,\ldots,N_{z}-2,N_{z}-1\right),
\]
where $h_{z}=z_{\max}^{\mathbf{III}}/N_{z}$ denotes the mesh-width.
Note that the point $z_{l_{1}l_{2}}$ is the center of the cell 
\[
R_{l_{1}l_{2}}=\left[l_{1}h_{z},\left(l_{1}+1\right)h_{z}\right]\times\left[l_{2}h_{z},\left(l_{2}+1\right)h_{z}\right].
\]
We then consider the approximation
\begin{equation}
\int_{\mathbb{R}_{*}^{2}}f\left(z,x,t\right)\ell\left(dz\right)\simeq\int_{R_{z}^{\mathbf{I}}}f\left(z,x,t\right)\ell\left(dz\right)+\int_{R_{z}^{\mathbf{II}}}f\left(z,x,t\right)\ell\left(dz\right)+\int_{R_{z}^{\mathbf{III}}}f\left(z,x,t\right)\ell\left(dz\right),\label{eq: integral disc - integral operator (1) truncated domain}
\end{equation}
where the individual terms on the right-hand side will be approximated
in different ways: the first one will be transformed into a diffusion
term by replacing the integrand function with its Taylor expansion;
for the second one, a particular quadrature formula is used that takes
into account the limiting singular behaviour of the L{\'e}vy measure as
$\left\Vert z\right\Vert \rightarrow0^{+}$; for the third one, a
generic method is used.

By substituting (\ref{eq: integral disc - Taylor approximation of the integrand})
in the first integral in (\ref{eq: integral disc - integral operator (1) truncated domain}),
it follows that
\begin{equation}
\int_{R_{z}^{\mathbf{I}}}f\left(z,x,t\right)\ell\left(dz\right)\simeq\frac{1}{2}\mathbf{1}^{\top}\left(u_{xx}\left(x,t\right)\circ I_{x}\left(\int_{R_{z}^{\mathbf{I}}}zz^{\top}\ell\left(dz\right)\right)I_{x}\right)\mathbf{1}.\label{eq: integral disc - integral operator (2) approximation over R0}
\end{equation}
Here, the entries of the matrix $\int_{R_{z}^{\mathbf{I}}}zz^{\top}\ell\left(dz\right)$
can be accurately approximated using a common numerical integrator.

Moving on to the second and third terms in (\ref{eq: integral disc - integral operator (1) truncated domain}),
we consider a quadrature formula of the form
\begin{equation}
\int_{R_{z}^{\mathbf{II}}\cup R_{z}^{\mathbf{III}}}f\left(z,x,t\right)\ell\left(dz\right)\simeq\sum_{l_{1},l_{2}=-N_{z}}^{N_{z}-1}\omega_{l_{1}l_{2}}f\left(z_{l_{1}l_{2}},x,t\right).\label{eq: integral disc - integral operator (3) approximation over R1 and R2}
\end{equation}
Defining the coefficients 
\begin{equation}
\omega_{l_{1}l_{2}}=\begin{cases}
0 & \text{if }l_{1},l_{2}:z_{l_{1}l_{2}}\in R_{z}^{\mathbf{I}},\\
\left\Vert z_{l_{1}l_{2}}\right\Vert ^{-2}\int_{R_{l_{1}l_{2}}}\left\Vert z\right\Vert ^{2}\ell\left(dz\right) & \text{if }l_{1},l_{2}:z_{l_{1}l_{2}}\in R_{z}^{\mathbf{II}},\\
\ell\left(z_{l_{1}l_{2}}\right)h_{z}^{2} & \text{if }l_{1},l_{2}:z_{l_{1}l_{2}}\in R_{z}^{\mathbf{III}},
\end{cases}\label{eq: integral disc - omega coefficients}
\end{equation}
a satisfactory level of accuracy is achieved, despite the integrand being
singular at the origin. Clearly, the quadrature weights used in $R_{z}^{\mathbf{II}}$
are constructed as integrals of the L{\'e}vy measure, which turns out
to be beneficial for the convergence behaviour (as $N_{z}\rightarrow\infty$).
Analogously to the entries of the matrix $\int_{R_{z}^{\mathbf{I}}}zz^{\top}\ell\left(dz\right)$
in (\ref{eq: integral disc - integral operator (2) approximation over R0}),
the integrals $\int_{R_{l_{1}l_{2}}}\left\Vert z\right\Vert ^{2}\ell\left(dz\right)$
can be precomputed using a common numerical integrator. Regarding
$R_{z}^{\mathbf{III}}$, the coefficients are obtained by applying
the classical midpoint rule, see for example \citet{quarteroni2007}.
Finally, note that the weights $\omega_{l_{1}l_{2}}$ are null over
$R_{z}^{\mathbf{I}}$, as the first integral in (\ref{eq: integral disc - integral operator (1) truncated domain})
has already been approximated through (\ref{eq: integral disc - integral operator (2) approximation over R0}).

For the above approximation of the integral term, we have the following convergence result.
\begin{prop}
\label{prop: integral disc - Quadrature scheme for the integral operator}
Let $u$ (resp.\ $\ell$) be sufficiently smooth with respect to $x$ (resp.\ $z$).
Let $z_{\max}^{\mathbf{I}}$ be directly proportional to $h_{z}$ and let $z_{\max}^{\mathbf{II}}$, $z_{\max}^{\mathbf{III}}$
be fixed (independent of $h_{z}$).
Let $\left\Vert \cdot\right\Vert $ be the Euclidean norm.
Assume $h_{z}$ is such that $z_{\max}^{\mathbf{II}}/h_{z}$ is an integer.
Then, for any given $\epsilon>0$ and $x_{\max}>0$, the quadrature error
\begin{align}
 E\left(x,t\right) & =\int_{\left\Vert z\right\Vert _{\infty}<z_{\max}^{\mathbf{III}}}f\left(z,x,t\right)\ell\left(dz\right)\label{eq: error formula}\\
 &\quad  -\frac{1}{2}\mathbf{1}^{\top}\left(u_{xx}\left(x,t\right)\circ I_{x}\left(\int_{R_{z}^{\mathbf{I}}}zz^{\top}\ell\left(dz\right)\right)I_{x}\right)\mathbf{1}-\sum_{l_{1},l_{2}=-N_{z}}^{N_{z}-1}\omega_{l_{1}l_{2}}f\left(z_{l_{1}l_{2}},x,t\right)\nonumber 
\end{align}
satisfies
\[
E\left(x,t\right)=\begin{cases}
O\left(h_{z}^{2-A_{\ell}}\right) & \text{if } 0<A_{\ell}<2,\\
O\left(h_{z}^{2-\epsilon}\right) & \text{if } A_{\ell}=0,
\end{cases}
\]
uniformly in $\left(x,t\right)\in \left[0,x_{\max}\right]^2 \times\left[0,T\right]$.
\end{prop}

\begin{proof}
 See Appendix \ref{sec: appendix - Modified midpoint rule}.
 \end{proof}

Proposition \ref{prop: integral disc - Quadrature scheme for the integral operator} corresponds to a result derived by \citet{wang2007} for the 1-dimensional setting.
Clearly, in our 2-dimensional setting, essentially second-order convergence is obtained for $A_{\ell} =0$ and the order of convergence equals $2-A_{\ell}$ whenever $0< A_{\ell} < 2$.
This forms a positive result. We note, however, that it is less advantageous compared to the 1-dimensional setting, as in this case \citet{wang2007} proved second-order convergence whenever $0\le A_{\ell} <1$, essentially second-order convergence if $A_{\ell}=1$, and an order of convergence equal to $3-A_{\ell}$ whenever $1<A_{\ell}<2$.
The less advantageous convergence result in the 2-dimensional setting is attributed to a lack of smoothness w.r.t.~$z$ of the scaled function $f\left(z,x,t\right)\left\Vert z\right\Vert ^{-2}$ in the neighbourhood of $z=0$, see the proof in Appendix \ref{sec: appendix - Modified midpoint rule}.

Using \eqref{eq: integral disc - integral operator (2) approximation over R0} and substituting (\ref{eq: intro - integrand function}) into (\ref{eq: integral disc - integral operator (3) approximation over R1 and R2}),
we define an approximating operator $\mathcal{A}_{\omega}$ and a number $r_{\omega}$ such that 
\[
\mathcal{A}u\left(x,t\right)-ru\left(x,t\right)\simeq\mathcal{A}_{\omega}u\left(x,t\right)-r_{\omega}u\left(x,t\right)\qquad\text{for any }\left(x,t\right)\in\mathbb{R}_{\geq0}^{2}\times\left(0,T\right],
\]
with
\begin{align}
\mathcal{A}_{\omega}u\left(x,t\right) & =\mu_{\omega}\left(x\right)^{\top}u_{x}\left(x,t\right)+\frac{1}{2}\mathbf{1}^{\top}\left(u_{xx}\left(x,t\right)\circ\Sigma_{\omega}\Sigma_{\omega}^{\top}\left(x\right)\right)\mathbf{1}+\left(\mathcal{B}_{\omega}u\right)\left(x,t\right),\label{eq: operator Aw - definition}\\
r_{\omega} & =r+\sum_{l_{1},l_{2}=-N_{z}}^{N_{z}-1}\omega_{l_{1}l_{2}},\nonumber 
\end{align}
where, for $i,j=1,2$,
\begin{align}
\mu_{\omega}^{\left(i\right)}\left(x\right) & =x^{\left(i\right)}\kappa_{\omega}^{\left(i\right)},\nonumber \\
\kappa_{\omega}^{\left(i\right)} & =r-\sum_{l_{1},l_{2}=-N_{z}}^{N_{z}-1}\omega_{l_{1}l_{2}}\left(e^{z_{l_{1}l_{2}}^{\left(i\right)}}-1\right),\nonumber \\
\left(\Sigma_{\omega}\Sigma_{\omega}^{\top}\right)^{\left(i,j\right)}\left(x\right) & =x^{\left(i\right)}x^{\left(j\right)}\left(\sigma_{\omega}\sigma_{\omega}^{\top}\right)^{\left(i,j\right)},\nonumber \\
\sigma_{\omega}\sigma_{\omega}^{\top} & =\sigma\sigma^{\top}+\int_{R_{z}^{\mathbf{I}}}zz^{\top}\ell\left(dz\right),\nonumber \\
\left(\mathcal{B}_{\omega}u\right)\left(x,t\right) & =\sum_{l_{1},l_{2}=-N_{z}}^{N_{z}-1}\omega_{l_{1}l_{2}}u\left(x+\gamma\left(z_{l_{1}l_{2}},x\right),t\right).\label{eq: integral disc - summation operator (1) definition}
\end{align}
Then, we approximate the solution $u$ of (\ref{eq: intro - IVP (1)})
by the function $v:\mathbb{R}_{\geq0}^{2}\times\left[0,T\right]\rightarrow\mathbb{R}$
which solves the following problem
\begin{equation}
\begin{cases}
v_{t}\left(x,t\right)=\mathcal{A}_{\omega}v\left(x,t\right)-r_{\omega}v\left(x,t\right) & \text{for any }\left(x,t\right)\in\mathbb{R}_{\geq0}^{2}\times\left(0,T\right]\\
v\left(x,0\right)=\phi\left(x\right).
\end{cases}\label{eq: integral disc - IVP (2)}
\end{equation}

\subsection{\label{subsec: spatial discretization}Spatial discretization}

In this section, we successively consider the spatial discretization
of the diffusion and summation terms in the operator $\mathcal{A}_{\omega}$.
The convection term will be discussed in Section \ref{sec: semi-Lagrangian theta-method}.

Let $R_{x}=\left[0,x_{\max}\right]\times\left[0,x_{\max}\right]$
be the truncated $x$-domain over which the solution to (\ref{eq: integral disc - IVP (2)})
is approximated and $N_{x}\in\mathbb{N}$ be a given number of spatial
grid points. Here, $x_{\max}$ is chosen heuristically in such a way
that the localization error is negligible. We construct a spatial
grid $\mathbf{x}$ in $R_{x}$ by applying, in each dimension, a strictly
increasing and smooth transformation $\varphi$ to an artificial uniform
grid. Let
\[
x_{m}=\varphi\left(\varphi^{-1}\left(0\right)+\frac{\varphi^{-1}\left(x_{\max}\right)-\varphi^{-1}\left(0\right)}{N_{x}}m\right)\qquad\left(m=0,1,\ldots,N_{x}\right).
\]
The elements of $\mathbf{x}$ are defined
by 
\[
x_{m_{1}m_{2}}=\left(x_{m_{1}},x_{m_{2}}\right)\qquad\left(m_{1},m_{2}=0,1,\ldots,N_{x}\right).
\]
The function $\varphi$ will be chosen in such a way that relatively
many points are placed in a region of financial and numerical interest. 

In what follows, we denote the values over $\mathbf{x}$ of any given
function $g:R_{x}\times\left[0,T\right]\rightarrow\mathbb{R}$ by
the vector
\begin{equation}
g\left(\mathbf{x},t\right)=\left[g\left(x_{00},t\right),g\left(x_{10},t\right),\ldots,g\left(x_{N_{x}-1,N_{x}},t\right),g\left(x_{N_{x}N_{x}},t\right)\right]^{\top}.\label{eq: spat, ordered spatial verse}
\end{equation}

\subsubsection{\label{subsec: spatial discretization - diffusion term}Diffusion
term}

In this subsection, we construct a semi-discrete diffusion matrix
$D$ such that
\begin{equation}
Dv\left(\mathbf{x},t\right)\simeq\left[\frac{1}{2}\mathbf{1}^{\top}\left(v_{xx}\left(x_{m_{1}m_{2}},t\right)\circ\Sigma_{\omega}\Sigma_{\omega}^{\top}\left(x_{m_{1}m_{2}}\right)\right)\mathbf{1}\right]_{m_{1},m_{2}=0,1,\ldots,N_{x}},\label{eq: spatial disc - matrix D (1)}
\end{equation}
where the right-hand side is a vector, whose elements are ordered
according to (\ref{eq: spat, ordered spatial verse}).

To this purpose, in each spatial dimension, we approximate the first-
and second-order derivatives of a given smooth function $g:\mathbb{R}\rightarrow\mathbb{R}$
by the following second-order central finite difference schemes 
\begin{align*}
g^{\prime}\left(x_{m}\right) & \simeq\alpha_{m}^{\left(-1\right)}g\left(x_{m-1}\right)+\alpha_{m}^{\left(0\right)}g\left(x_{m}\right)+\alpha_{m}^{\left(1\right)}g\left(x_{m+1}\right)\\
g^{\prime\prime}\left(x_{m}\right) & \simeq\beta_{m}^{\left(-1\right)}g\left(x_{m-1}\right)+\beta_{m}^{\left(0\right)}g\left(x_{m}\right)+\beta_{m}^{\left(1\right)}g\left(x_{m+1}\right)
\end{align*}
with coefficients
\[
\alpha_{m}^{\left(-1\right)}=\frac{-h_{x,m+1}}{h_{x,m}\left(h_{x,m}+h_{x,m+1}\right)},\qquad\alpha_{m}^{\left(0\right)}=\frac{h_{x,m+1}-h_{x,m}}{h_{x,m}h_{x,m+1}},\qquad\alpha_{m}^{\left(1\right)}=\frac{h_{x,m}}{h_{x,m+1}\left(h_{x,m}+h_{x,m+1}\right)},
\]
\[
\beta_{m}^{\left(-1\right)}=\frac{2}{h_{x,m}\left(h_{x,m}+h_{x,m+1}\right)},\qquad\beta_{m}^{\left(0\right)}=\frac{-2}{h_{x,m}h_{x,m+1}},\qquad\beta_{m}^{\left(1\right)}=\frac{2}{h_{x,m+1}\left(h_{x,m}+h_{x,m+1}\right)},
\] where $h_{x,m}=x_{m}-x_{m-1}$. Concerning the boundary of the truncated spatial domain, we modify
the previous formulae in the following way. At the lower boundary
$x_{0}=0$, the first- and second-order derivative terms in (\ref{eq: integral disc - IVP (2)})
vanish. Hence, it is natural to choose $\alpha_{0}^{\left(j\right)}=0$
and $\beta_{0}^{\left(j\right)}=0$ for any $j=\left\{ -1,0,1\right\} $.
At the upper boundary $x_{N_{x}}=x_{\max}$, we make the natural assumption
that the solution $v$ behaves linearly in $x$, thus we choose $\beta_{N_{x}}^{\left(j\right)}=0$
for any $j=\left\{ -1,0,1\right\} $, and we approximate the first-order
derivative by the first-order backward finite difference scheme. 

Noting that $\mathbf{x}$ is the Cartesian product of two identical
1-dimensional grids, by employing the 1-directional finite difference
formulae in both the spatial dimensions, it leads to the matrix $D$
defined by 
\begin{equation}
D=\frac{1}{2}\left(\sigma_{\omega}\sigma_{\omega}^{\top}\right)^{\left(1,1\right)}I\otimes I_{\mathbf{x}}^{2}D_{2}+\left(\sigma_{\omega}\sigma_{\omega}^{\top}\right)^{\left(1,2\right)}I_{\mathbf{x}}D_{1}\otimes I_{\mathbf{x}}D_{1}+\frac{1}{2}\left(\sigma_{\omega}\sigma_{\omega}^{\top}\right)^{\left(2,2\right)}I_{\mathbf{x}}^{2}D_{2}\otimes I.\label{eq: spatial disc - matrix D (2)}
\end{equation}
Here, $I\in\mathbb{R}^{\left(N_{x}+1\right)\times\left(N_{x}+1\right)}$
is the identity matrix, $I_{\mathbf{x}}=\text{diag}\left(x_{0}^{\left(i\right)},\ldots,x_{N_{x}}^{\left(i\right)}\right)$
and $\otimes$ denotes the Kronecker product.\footnote{In this paper, we use the convention $AB\otimes CD=\left(AB\right)\otimes\left(CD\right)$,
for any suitable matrices $A,B,C,D$.} The matrices $D_{1},D_{2}\in\mathbb{R}^{\left(N_{x}+1\right)\times\left(N_{x}+1\right)}$
are the matrices representing numerical differentiation of first-
and second-order by the relevant finite difference formulae above.
The mixed derivative has been approximated by applying the finite
difference formula for the first-order derivative subsequently in
the two spatial dimensions.

\subsubsection{\label{subsec: spatial disc - summation term}Summation term}

In this section, we derive an efficient method to approximate the
summation term $\left(\mathcal{B}_{\omega}v\right)\left(\mathbf{x},t\right)$
given the values of $v\left(\mathbf{x},t\right)$. Unlike the differential
component of $\mathcal{A}_{\omega}$, we do not construct a matrix
$B_{\omega}$ such that $\left(\mathcal{B}_{\omega}v\right)\left(\mathbf{x},t\right)\simeq B_{\omega}v\left(\mathbf{x},t\right)$,
as this matrix would be large and dense.

Assuming that the values of $v$ are known for all $\left(x,t\right)\in R_{x}\times\left[0,T\right]$,
using formula (\ref{eq: integral disc - summation operator (1) definition})
to directly evaluate $\left(\mathcal{B}_{\omega}v\right)\left(\mathbf{x},t\right)$
would require $O\left(N_{x}^{2}N_{z}^{2}\right)$ elementary operations,
which is computationally too expensive. For this reason, a particularly
efficient method combining interpolation and FFT is considered, which
extends the approach by \citet{wang2007}.

Let $N_{y}^{-},N_{y}^{+}\in\mathbb{N}$ be any given natural numbers
and let $\mathbf{y}^{\text{out}}$ and $\mathbf{y}^{\text{in}}$ be
two grids of points defined~by\footnote{The superscripts stand for ``input'' and ``output''.}
\begin{align*}
y_{m_{1}m_{2}}^{\text{out}} & =\left(e^{m_{1}h_{z}},e^{m_{2}h_{z}}\right)\qquad\left(m_{1},m_{2}=-N_{y}^{-},-N_{y}^{-}+1,\ldots,N_{y}^{+}-1,N_{y}^{+}\right),\\
y_{m_{1}m_{2}}^{\text{in}} & =\left(e^{\left(m_{1}+\frac{1}{2}\right)h_{z}},e^{\left(m_{2}+\frac{1}{2}\right)h_{z}}\right)\qquad\left(m_{1},m_{2}=-N_{z}-N_{y}^{-},-N_{z}-N_{y}^{-}+1,\ldots,N_{z}+N_{y}^{+}-2,N_{z}+N_{y}^{+}-1\right),
\end{align*}
then it holds that
\begin{equation}
\left(\mathcal{B}_{\omega}v\right)\left(y_{m_{1}m_{2}}^{\text{out}},t\right)=\sum_{l_{1},l_{2}=-N_{z}}^{N_{z}-1}\omega_{l_{1}l_{2}}v\left(y_{l_{1}+m_{1},l_{2}+m_{2}}^{\text{in}},t\right)\qquad\left(m_{1},m_{2}=-N_{y}^{-},-N_{y}^{-}+1,\ldots,N_{y}^{+}-1,N_{y}^{+}\right).\label{eq: spatial disc - summation operator (2) over the y-grids}
\end{equation}
Clearly, the summation term \eqref{eq: spatial disc - summation operator (2) over the y-grids} can be viewed as a discrete 2-dimensional cross-correlation. It is well known, see for instance \citet[Chapter 3]{plonka2018},  that it can be written in the form
\begin{equation}
\left(\mathcal{B}_{\omega}v\right)\left(\mathbf{y}^{\text{out}},t\right)=\tilde{I}Cv\left(\mathbf{y}^{\text{in}},t\right)\label{eq: spatial disc - summation operator (3) over the y-grids}
\end{equation}
where:
\begin{itemize}
\item $C\in\mathbb{R}^{\left(\sharp\text{in}\right)^{2}\times\left(\sharp\text{in}\right)^{2}}$
is a circulant matrix whose first row is given by $C_{1,\cdot}^{\top}$
with 
\begin{equation}
C_{1,\cdot}=\text{vec}\left(\left[\begin{array}{cc}
\Omega & 0_{\sharp\mathbf{z}\times\left(\sharp\text{in}-\sharp\mathbf{z}\right)}\\
0_{\left(\sharp\text{in}-\sharp\mathbf{z}\right)\times\sharp\mathbf{z}} & 0_{\left(\sharp\text{in}-\sharp\mathbf{z}\right)\times\left(\sharp\text{in}-\sharp\mathbf{z}\right)}
\end{array}\right]\right).\label{eq: spatial disc - matrix C (1) first row}
\end{equation}
Here, $0_{P\times Q}$ denotes the null matrix of dimensions $P\times Q$,
$\text{vec}\left(\cdot\right)$ denotes the vectorization of a matrix,
$\sharp$ indicates the number of points of a given grid in one direction
and $\Omega\in\mathbb{R}^{\sharp\mathbf{z}\times\sharp\mathbf{z}}$
is the matrix whose entries are the coefficients $\omega_{l_{1}l_{2}}$
defined by (\ref{eq: integral disc - omega coefficients}). For an example of a matrix $C$, we refer to Appendix \ref{sec: appendix - summation operator as a circulant matrix-vector multiplication}.

The matrix-vector
multiplication $Ca$, for any given vector $a\in\mathbb{R}^{\left(\sharp\text{in}\right)^{2}\times1}$,
can be obtained by two (1-dimensional) FFTs and one (1-dimensional)
inverse FFT, requiring just $O\left(\left(\sharp\text{in}\right)^{2}\cdot\log\sharp\text{in}\right)$
elementary operations. The pertinent formula is
\begin{equation}
Ca=\text{ifft}\left(\text{fft}\left(C_{1,\cdot}\right)^{H}\circ\text{fft}\left(a\right)\right),\label{eq: spatial disc - matrix C (2) FFT}
\end{equation}
where $^{H}$ denotes the complex conjugate. 
\item $\tilde{I}\in\mathbb{R}^{\left(\sharp\text{out}\right)^{2}\times\left(\sharp\text{in}\right)^{2}}$
is obtained from the identity matrix $I\in\mathbb{R}^{\left(\sharp\text{in}\right)^{2}\times\left(\sharp\text{in}\right)^{2}}$
by removing the rows corresponding to the zeros in the following vector
\[
\text{vec}\left(\left[\begin{array}{cc}
1_{\sharp\text{out}\times\sharp\text{out}} & 0_{\sharp\text{out}\times\left(\sharp\text{in}-\sharp\text{out}\right)}\\
0_{\left(\sharp\text{in}-\sharp\text{out}\right)\times\sharp\text{out}} & 0_{\left(\sharp\text{in}-\sharp\text{out}\right)\times\left(\sharp\text{in}-\sharp\text{out}\right)}
\end{array}\right]\right).
\]
Here, $1_{P\times P}$ denotes a $P\times P$ matrix whose elements
are all equal to 1. We note that the matrix-vector multiplication
$Cv\left(\mathbf{y}^{\text{in}},t\right)$ in (\ref{eq: spatial disc - summation operator (3) over the y-grids})
returns a value also for grid points that can be discarded. The purpose
of $\tilde{I}$ is precisely to extract only those entries that correspond
to $\left(\mathcal{B}_{\omega}v\right)\left(\mathbf{y}^{\text{out}},t\right)$.
\end{itemize}
In order to obtain an approximation to $\left(\mathcal{B}_{\omega}v\right)\left(\mathbf{x},t\right)$
using (\ref{eq: spatial disc - summation operator (3) over the y-grids}),
we need to interpolate both the input and the output value in (\ref{eq: spatial disc - summation operator (3) over the y-grids})
since $\mathbf{y}^{\text{in}}$ and $\mathbf{y}^{\text{out}}$ are
generally different from $\mathbf{x}$. Let $T^{\text{in}}\in\mathbb{R}^{\left(\sharp\text{in}\right)^{2}\times\left(N_{x}+1\right)^{2}}$
be a matrix representing an interpolation procedure from the $\mathbf{x}$
grid to the $\mathbf{y}^{\text{in}}$ grid and let $T^{\text{out}}\in\mathbb{R}^{\left(N_{x}+1\right)^{2}\times\left(\sharp\text{in}\right)^{2}}$
be a matrix representing an interpolation procedure from the $\mathbf{y}^{\text{out}}$
grid to the $\mathbf{x}$ grid. Then
\begin{align}
v\left(\mathbf{y}^{\text{in}},t\right) & \simeq T^{\text{in}}v\left(\mathbf{x},t\right),\label{eq: spatial disc - interpolation input value}\\
\left(\mathcal{B}_{\omega}v\right)\left(\mathbf{x},t\right) & \simeq T^{\text{out}}\left(\mathcal{B}_{\omega}v\right)\left(\mathbf{y}^{\text{out}},t\right).\label{eq: spatial disc - interpolation output value}
\end{align}
Note that, by using Lagrange interpolation, the interpolation matrices
are sparse and have at most $P+1$ nonzero entries per row, where
$P$ is the polynomial degree. Let $M$ be the number of rows, it
follows that the corresponding matrix--vector multiplications require
a number of operations of order $O\left(MP\right)$, and are therefore
negligible compared with multiplication performed via FFT.

From (\ref{eq: spatial disc - summation operator (3) over the y-grids}),
(\ref{eq: spatial disc - interpolation input value}) and (\ref{eq: spatial disc - interpolation output value}),
we arrive at the approximation 
\begin{equation}
\left(\mathcal{B}_{\omega}v\right)\left(\mathbf{x},t\right)\simeq B_{\omega}v\left(\mathbf{x},t\right),\label{eq: spatial disc - Matrix Bw}
\end{equation}
where $B_{\omega}\in\mathbb{R}^{\left(N_{x}+1\right)^{2}\times\left(N_{x}+1\right)^{2}}$
is given by 
\begin{equation}
B_{\omega}=T^{\text{out}}\tilde{I}CT^{\text{in}}.\label{eq: spatial disc - Matrix Bw 2}
\end{equation}
We emphasize that $B_{\omega}$ is only used for notational purposes
and never explicitly computed. To compute the right-hand side of (\ref{eq: spatial disc - Matrix Bw}),
we always use
\begin{equation}
B_{\omega}V=T^{\text{out}}\tilde{I}\,\text{ifft}\left(\text{fft}\left(C_{1,\cdot}\right)^{H}\circ\text{fft}\left(T^{\text{in}}V\right)\right),\label{eq: spatial disc - matrix S (2) multiplication by FFT}
\end{equation}
for any vector $V\in\mathbb{R}^{\left(N_{x}+1\right)^{2}\times1}$.
Figure 
\begin{figure}[t]
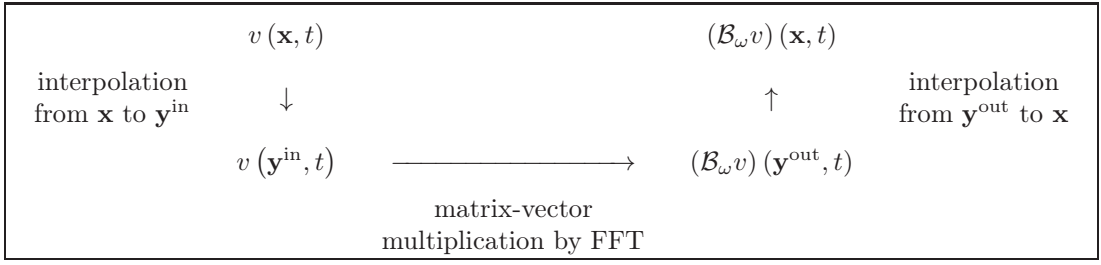

\caption{\label{fig: FFT-interpolation procedure to compute explicitly the cyclic summation}Diagram
of the scheme used to approximate $\left(\mathcal{B}_{\omega}v\right)\left(\mathbf{x},t\right)$}

\centering{}\medskip{}
\begin{tabular}{|ccccc|}
\hline 
$\begin{array}{c}
\\
\\
\end{array}$ & $v\left(\mathbf{x},t\right)$ &  & $\left(\mathcal{B}_{\omega}v\right)\left(\mathbf{x},t\right)$ & \tabularnewline
$\begin{array}{c}
\text{interpolation}\\
\text{from \ensuremath{\mathbf{x}} to }\mathbf{y}^{\text{in}}
\end{array}$ & $\downarrow$ &  & $\uparrow$ & $\begin{array}{c}
\text{interpolation}\\
\text{from \ensuremath{\mathbf{y}^{\text{out}}} to }\mathbf{x}
\end{array}$\tabularnewline
$\begin{array}{c}
\\
\\
\end{array}$ & $v\left(\mathbf{y}^{\text{in}},t\right)$ & $\xrightarrow{\text{\hspace{3cm}}}$ & $\left(\mathcal{B}_{\omega}v\right)\left(\mathbf{y}^{\text{out}},t\right)$ & \tabularnewline
 &  & $\begin{array}{c}
\text{matrix-vector}\\
\text{multiplication by FFT}
\end{array}$ &  & \tabularnewline
\hline 
\end{tabular}
\end{figure}
\ref{fig: FFT-interpolation procedure to compute explicitly the cyclic summation}
provides a schematic illustration of how FFT and interpolation are
combined to evaluate (\ref{eq: spatial disc - matrix S (2) multiplication by FFT}).

\subsubsection{\label{subsec: Cell averaging}Cell averaging}

We conclude the spatial discretization with a technique for handling
the non-smoothness of the initial function $\phi$ of (\ref{eq: intro - IVP (1)}).
As it turns out, pointwise valuation of $\phi$ over the spatial
grid can lead to deteriorated (spatial) convergence behaviour, which
can be alleviated by applying cell averaging.

Let
\begin{align*}
x_{m+\frac{1}{2}} & =\frac{1}{2}\left(x_{m}+x_{m+1}\right)\qquad\left(m=0,1,\ldots,N_{x}-1\right)\\
h_{x,m+\frac{1}{2}} & =x_{m+\frac{1}{2}}-x_{m-\frac{1}{2}}\qquad\left(m=0,1,\ldots,N_{x}\right)
\end{align*}
with $x_{-\frac{1}{2}}=-x_{\frac{1}{2}}$ and $x_{N_{x}+\frac{1}{2}}=2x_{\max}-x_{N_{x}-\frac{1}{2}}$.
Then, we use the approximation
\begin{equation}
v\left(x_{m_{1}m_{2}},0\right)\simeq\frac{1}{h_{x,m_{1}+\frac{1}{2}}h_{x,m_{2}+\frac{1}{2}}}\int_{x_{m_{1}-\frac{1}{2}}}^{x_{m_{1}+\frac{1}{2}}}\int_{x_{m_{2}-\frac{1}{2}}}^{x_{m_{2}+\frac{1}{2}}}\phi\left(x_{1},x_{2}\right)dx_{2}dx_{1},\label{eq: spatial disc - cell-averaging}
\end{equation}
whenever the cell $\left[x_{m_{1}-\frac{1}{2}},x_{m_{1}+\frac{1}{2}}\right)\times\left[x_{m_{2}-\frac{1}{2}},x_{m_{2}+\frac{1}{2}}\right)$
has a nonempty intersection with the set of non-smoothness of $\phi$.

\subsection{\label{sec: semi-Lagrangian theta-method}Temporal discretization:
the semi-Lagrangian $\theta$-method}

The problem (\ref{eq: integral disc - IVP (2)}) can be convection-dominated.
To account for this, we shall consider temporal discretization using
the $\theta$-method combined with the semi-Lagrangian approach, as
described by \citet{spiegelman2006}. The semi-Lagrangian method follows,
in each time step, the characteristics backwards in time to determine
the departure points of the spatial grid points.

Let $x:\left[0,T\right]\rightarrow\mathbb{R}_{\geq0}^{2}$ and $v^{*}:\left[0,T\right]\rightarrow\mathbb{R}$
such that $v^{*}\left(t\right)=v\left(x\left(t\right),t\right)$.
The derivative of $v^{*}$ is given by
\[
v_{t}^{*}\left(t\right)=v_{t}\left(x\left(t\right),t\right)+x_{t}\left(t\right)^{\top}v_{x}\left(x\left(t\right),t\right).
\]
Assume $x$ satisfies the following (linear) ODE: 
\begin{equation}
x_{t}\left(t\right)=-\mu_{\omega}\left(x\left(t\right)\right)\qquad\left(0<t\leq T\right).\label{eq: temporal disc - SL trajectory (1) ODE}
\end{equation}
Then
\begin{equation}
v_{t}^{*}\left(t\right)=\left(\mathcal{A}_{\omega}^{\text{SL}}-r_{\omega}\right)v\left(x\left(t\right),t\right)\qquad\left(0<t\leq T\right),\label{eq: temporal disc - IVP (3)}
\end{equation}
where 
\[
\mathcal{A}_{\omega}^{\text{SL}}v\left(x,t\right)=\frac{1}{2}\mathbf{1}^{\top}\left(v_{xx}\left(x,t\right)\circ\Sigma_{\omega}\Sigma_{\omega}^{\top}\left(x\right)\right)\mathbf{1}+\left(\mathcal{B}_{\omega}v\right)\left(x,t\right).
\]
Clearly, $\mathcal{A}_{\omega}^{\text{SL}}$ is obtained from $\mathcal{A}_{\omega}$
by omitting the convection term. 

Let parameter $\theta\in\left[0,1\right]$. Let $\mathbf{t}=\left(t_{n}\right)_{n=0}^{N_{t}}$
be any given uniform grid with step size $h_{t}=\frac{T}{N_{t}}$.
For any given $n=1,2,\ldots,N_{t}$, approximating (\ref{eq: temporal disc - IVP (3)})
using the $\theta$-method and substituting the definition of $v^{*}$,
we obtain
\begin{equation}
\frac{v\left(x\left(t_{n}\right),t_{n}\right)-v\left(x\left(t_{n-1}\right),t_{n-1}\right)}{h_{t}}\simeq\theta\left(\mathcal{A}_{\omega}^{\text{SL}}-r_{\omega}\right)v\left(x\left(t_{n}\right),t_{n}\right)+\left(1-\theta\right)\left(\mathcal{A}_{\omega}^{\text{SL}}-r_{\omega}\right)v\left(x\left(t_{n-1}\right),t_{n-1}\right).\label{eq: temporal disc - mid formula}
\end{equation}

The approximation (\ref{eq: temporal disc - mid formula}) holds along
any trajectory satisfying (\ref{eq: temporal disc - SL trajectory (1) ODE}).
In each given time step from $t_{n-1}$ to $t_{n}$, the semi-Lagrangian
approach involves selecting the set of trajectories that intersect
the points $\left(\mathbf{x},t_{n}\right)$, ensuring that an approximation
is defined on the fixed grid $\mathbf{x}$. Let $\mathbf{x}^{\text{SL}}$
denote the grid corresponding to $t_{n-1}$ along this set of trajectories.
Its elements are given by $x_{m_{1}m_{2}}^{\text{SL}}=\left(x_{m_{1}}^{\text{SL}},x_{m_{2}}^{\text{SL}}\right)$
where $x_{m_{i}}^{\text{SL}}$ is obtained by (\ref{eq: temporal disc - SL trajectory (1) ODE})
as 
\[
x_{m_{i}}^{\text{SL}}=x_{m_{i}}e^{\kappa_{\omega}^{\left(i\right)}h_{t}}\qquad\left(m_{i}=0,1,\ldots,N_{x}\right).
\]
Then (\ref{eq: temporal disc - mid formula}) becomes 
\[
\frac{v\left(\mathbf{x},t_{n}\right)-v\left(\mathbf{x}^{\text{SL}},t_{n-1}\right)}{h_{t}}\simeq\theta\left(\mathcal{A}_{\omega}^{\text{SL}}-r_{\omega}\right)v\left(\mathbf{x},t_{n}\right)+\left(1-\theta\right)\left(\mathcal{A}_{\omega}^{\text{SL}}-r_{\omega}\right)v\left(\mathbf{x}^{\text{SL}},t_{n-1}\right).
\]
Interpolation is employed to acquire approximations at the grid $\mathbf{x}^{\text{SL}}$.
Let $T^{\text{SL}}\in\mathbb{R}^{\left(N_{x}+1\right)^{2}\times\left(N_{x}+1\right)^{2}}$
be a matrix representing an interpolation procedure from the $\mathbf{x}$
grid to the $\mathbf{x}^{\text{SL}}$ grid. Together with the approximation
of the diffusion and summation terms, discussed in Section \ref{subsec: spatial discretization},
we obtain 
\begin{align}
v\left(\mathbf{x}^{\text{SL}},t_{n-1}\right) & \simeq T^{\text{SL}}v\left(\mathbf{x},t_{n-1}\right),\label{eq: temp, interpolation SL}\\
\left(\mathcal{A}_{\omega}^{\text{SL}}-r_{\omega}\right)v\left(\mathbf{x},t_{n}\right) & \simeq\left(D+B_{\omega}-r_{\omega}I\right)v\left(\mathbf{x},t_{n}\right),\nonumber \\
\left(\mathcal{A}_{\omega}^{\text{SL}}-r_{\omega}\right)v\left(\mathbf{x}^{\text{SL}},t_{n-1}\right) & \simeq T^{\text{SL}}\left(D+B_{\omega}-r_{\omega}I\right)v\left(\mathbf{x},t_{n-1}\right).\nonumber 
\end{align}

This leads to the following natural definition of the approximation
$V^{n}$ to the exact solution vector $v\left(\mathbf{x},t_{n}\right)$:
\begin{equation}
\left(I-h_{t}\theta\left(D+B_{\omega}-r_{\omega}I\right)\right)V^{n}=T^{\text{SL}}\left(I+h_{t}\left(1-\theta\right)\left(D+B_{\omega}-r_{\omega}I\right)\right)V^{n-1}\label{eq: temporal disc - IVP (4)}
\end{equation}
for $n=1,2,\ldots,N_{t}$. The initial vector $V^{0}$ is defined
by pointwise valuation on the spatial grid $\mathbf{x}$ of the pay-off
function $\phi$, except near the set of non-smoothness, where cell
averaging is employed (see Section \ref{subsec: Cell averaging}).
The time-stepping scheme (\ref{eq: temporal disc - IVP (4)}) is called
the semi-Lagrangian $\theta$-method. We shall apply (\ref{eq: temporal disc - IVP (4)})
with $\theta=\frac{1}{2}$, which is also called the semi-Lagrangian
Crank--Nicolson method. Here, to account for the non-smoothness of
$\phi$, a damping procedure is used where the first time step (i.e.
$n=1$) is replaced by four time steps of size equal to $\frac{1}{4}h_{t}$
of (\ref{eq: temporal disc - IVP (4)}) with $\theta=1$.

It remains to consider the treatment of the discretized integral term
in (\ref{eq: temporal disc - IVP (4)}), represented formally by the
matrix $B_{\omega}$. Recall from Section \ref{subsec: spatial disc - summation term}
that $B_{\omega}$ is never actually computed. To effectively handle
this term, we shall employ fixed-point iteration:
\begin{equation}
\left(I-h_{t}\theta\left(D-r_{\omega}I\right)\right)Y^{n,k}=h_{t}\theta B_{\omega}Y^{n,k-1}+T^{\text{SL}}\left(I+h_{t}\left(1-\theta\right)\left(D-r_{\omega}I\right)\right)V^{n-1}+h_{t}\left(1-\theta\right)T^{\text{SL}}B_{\omega}V^{n-1}\label{eq: temporal disc - fixed-point (1) recursive equation}
\end{equation}
for $k=1,2,\ldots$. Here matrix-vector multiplications involving
$B_{\omega}$ are always computed by the efficient FFT algorithm of
Section \ref{subsec: spatial disc - summation term}. For a given
tolerance $tol>0$ sufficiently small, we use the following stopping
criterion
\begin{equation}
\max_{m_{1},m_{2}}\frac{\left|Y_{m_{1}m_{2}}^{n,k}-Y_{m_{1}m_{2}}^{n,k-1}\right|}{\max\left\{ 1,\left|Y_{m_{1}m_{2}}^{n,k}\right|\right\} }<tol\label{eq: temporal disc - fixed-point (3) exit condition}
\end{equation}
 and define $V^{n}=Y^{n,k}$. 

The starting vector $Y^{n,0}$ for the fixed-point iteration is commonly
chosen in the literature as $Y^{n,0}=V^{n-1}$. Here, we shall consider
a more accurate starting vector, defined by higher-order extrapolation
from known approximations at previous temporal grid points:
\begin{equation}
Y^{n,0}=\begin{cases}
V^{n-1} & \textrm{if } n=1,\\
2V^{n-1}-V^{n-2} & \textrm{if } n=2,\\
3V^{n-1}-3V^{n-2}+V^{n-3} & \textrm{if } n=3,\\
4V^{n-1}-6V^{n-2}+4V^{n-3}-V^{n-4} & \textrm{if } n\geq4.
\end{cases}\label{eq: temporal disc - fixed-point (2) extrapolation starting value}
\end{equation}
This yields a significant reduction in the number of fixed-point iterations
compared to the common choice.

Finally, for the linear system in (\ref{eq: temporal disc - fixed-point (1) recursive equation})
we apply the BiCGSTAB iterative solver using an ILU preconditioner.

Our complete algorithm for the numerical solution of problem (\ref{eq: intro - IVP (1)})
is outlined in Algorithm \ref{alg: template of the algorithm}.
\begin{algorithm}[t]
\caption{\label{alg: template of the algorithm}Outline of the algorithm}

\textbf{precomputations:}
\begin{description}
\item [{$\bullet$}] define the grids $\mathbf{z}$, $\mathbf{x}$, $\mathbf{y}^{\text{in}}$,
$\mathbf{y}^{\text{out}}$, $\mathbf{t}$ and $\mathbf{x}^{\text{SL}}$
\item [{$\bullet$}] define the matrix $D$ given by (\ref{eq: spatial disc - matrix D (2)})
and compute the ILU factorization of $I-h_{t}\theta\left(D-r_{\omega}I\right)$
\item [{$\bullet$}] define the vector $C_{1,\cdot}$ given by (\ref{eq: spatial disc - matrix C (1) first row})
and compute $\text{fft}\left(C_{1,\cdot}\right)$
\item [{$\bullet$}] define the matrices $T^{\text{in}}$, $T^{\text{out}}$
and $T^{\text{SL}}$ given by (\ref{eq: spatial disc - interpolation input value}),
(\ref{eq: spatial disc - interpolation output value}) and (\ref{eq: temp, interpolation SL})
\item [{$\bullet$}] choose $\theta=\frac{1}{2}$
\end{description}
\textbf{time-stepping:}
\begin{description}
\item [{$\phantom{\bullet}$}] compute $V^{0}=\phi\left(\mathbf{x}\right)$
and apply cell averaging (\ref{eq: spatial disc - cell-averaging})
\item [{$\phantom{\bullet}$}] for $n=1,2,\ldots,N_{t}$
\begin{description}
\item [{1.}] compute $B_{\omega}V^{n-1}$ using (\ref{eq: spatial disc - matrix S (2) multiplication by FFT})
\item [{2.}] compute $W^{n-1}=T^{\text{SL}}\left(I+h_{t}\left(1-\theta\right)\left(D-r_{\omega}I\right)\right)V^{n-1}+h_{t}\left(1-\theta\right)T^{\text{SL}}B_{\omega}V^{n-1}$
\item [{3.}] compute $Y^{n,0}$ given by (\ref{eq: temporal disc - fixed-point (2) extrapolation starting value})
\item [{4.}] for $k=1,2,\ldots$
\begin{description}
\item [{i.}] compute $B_{\omega}Y^{n,k-1}$ using (\ref{eq: spatial disc - matrix S (2) multiplication by FFT})
\item [{ii.}] solve $\left(I-h_{t}\theta\left(D-r_{\omega}I\right)\right)Y^{n,k}=h_{t}\theta B_{\omega}Y^{n,k-1}+W^{n-1}$
using BiCGSTAB
\end{description}
\item [{5.}] end for if $Y^{n,k}$ satisfies (\ref{eq: temporal disc - fixed-point (3) exit condition})
\item [{6.}] let $V^{n}=Y^{n,k}$
\end{description}
\item [{$\phantom{\bullet}$}] end for
\end{description}
\end{algorithm}

\section{\label{sec: Numerical experiment}Numerical experiments}

We consider an European put-on-the-average option, which has the pay-off
function
\[
\phi\left(x\right)=\max\left(K-\frac{1}{2}\left(x^{\left(1\right)}+x^{\left(2\right)}\right),0\right)
\]
with fixed strike price $K>0$. Clearly, $\phi$ is non-smooth over
the set $\left\{ x\in\mathbb{R}_{\geq0}^{2}:x^{\left(1\right)}+x^{\left(2\right)}=2K\right\} $.
To define the non-uniform grid $\mathbf{x}$, we use the same transformation $\varphi$ as in \citet{inthout2023}. Let $c$, $x_{\rm int}$ be two given positive numbers. We choose the function $\varphi$ in  Section \ref{subsec: spatial discretization}  as \[
\varphi\left(\xi\right)=\begin{cases}
c\xi & \textrm{if } 0\leq\xi\leq\xi_{\rm int},\\
x_{\rm int}+c\sinh\left(\xi-\xi_{\rm int}\right) & \textrm{if }  \xi_{\rm int}<\xi\leq\xi_{\rm max},
\end{cases}
\]
with 
\[
\xi_{\rm int}=\frac{x_{\rm int}}{c},\quad\xi_{\rm max}=\xi_{\rm int}+\sinh^{-1}\left(\frac{x_{\rm max}-x_{\rm int}}{c}\right).
\]
In this way, the resulting spatial grid in each direction is uniform over $\left[0,x_{\rm int}\right]$, whereas in $\left[x_{\rm int},x_{\rm max}\right]$ the distance between consecutive grid points gradually increases as one moves away from $x_{\rm int}$. The limit of the fraction of spatial grid points within the interval $\left[0,x_{\rm int}\right]$
as $N_{x}\rightarrow\infty$, denoted by $F$, is given
by 
\[
F =\frac{\xi_{\rm int}}{\xi_{\rm max}}=\left(1+\frac{c}{x_{\rm int}}\sinh^{-1}\left(\frac{x_{\rm max}-x_{\rm int}}{c}\right)\right)^{-1}.
\]
Note that $F\rightarrow\frac{x_{\rm int}}{x_{\rm max}}$
as $c\rightarrow\infty$, which corresponds to the uniform
case.

Moving on to the L{\'e}vy measure, we model the jump component in (\ref{eq: intro - asset prices dynamics})
by a pure-jump 2-dimensional Normal Tempered Stable process. It is
characterized by the parameters $0 \leq \alpha<1$, $\delta>0$, $\lambda>0$,
$\eta\in\mathbb{R}^{2\times1}$ and a positive definite symmetric
matrix $\rho\in\mathbb{R}^{2\times2}$. The case where $\alpha=0$ is known as Variance Gamma, while the case where $\alpha=\frac{1}{2}$ is known as Normal Inverse Gaussian. Both are commonly used to model financial dynamics. The L{\'e}vy measure is given
by
\[
\ell\left(z\right)=\frac{\delta}{\pi}\sqrt{\frac{\left(\left\Vert \eta\right\Vert _{\rho}^{2}+2\lambda\right)^{1+\alpha}}{\det\left[\rho\right]}}K_{1+\alpha}\left(\sqrt{\left\Vert \eta\right\Vert _{\rho}^{2}+2\lambda}\left\Vert z\right\Vert _{\rho}\right)\left\Vert z\right\Vert _{\rho}^{-1-\alpha}e^{\left\langle \eta,z\right\rangle _{\rho}}
\]
where $K_{\nu}\left( \tau \right)=\frac{1}{2}\int_{0}^{\infty}y^{\nu-1}e^{-\frac{1}{2}\tau \left(y+y^{-1}\right)}dy$, for $\tau>0$,  
denotes the modified Bessel function of the second kind\footnote{See \citet[Appendix A]{schoutens2003}.}, 
$\left\langle x,y\right\rangle _{\rho}=x^{\top}\rho^{-1}y$ and $\left\Vert x\right\Vert _{\rho}=\sqrt{\left\langle x,x\right\rangle _{\rho}}$
is its induced norm. The constants $A_{\ell},B_{\ell}$, $C_{\ell}$ and $C^{\prime}_{\ell}$
in (\ref{eq: intro - upper bound of the L=0000E9vy measure}) are
defined, with respect to $\left\Vert \cdot\right\Vert _{\rho}$, as
\begin{align*}
&A_{\ell}  =2\alpha, \qquad
B_{\ell}  =\sqrt{\left\Vert \eta\right\Vert _{\rho}^{2}+2\lambda}-\left\Vert \eta\right\Vert _{\rho},\\
&C_{\ell}\left(h\right)  =\delta\frac{2^{\alpha}\Gamma\left(\alpha+1\right)}{\pi\sqrt{\det\left[\rho\right]}}e^{h\left\Vert \eta\right\Vert _{\rho}},\\
&C_{\ell}^{\prime}\left(h\right)  =C_{\ell}\left(h\right)\left\Vert \rho^{\frac{1}{2}}\right\Vert _{\rho}\left\Vert \rho^{-1}\right\Vert _{\rho}\left(h\left\Vert \eta\right\Vert _{\rho}+\left(A_{\ell}+2\right)\right).
\end{align*} The variance of the random variable $L\left(t\right)=\int_{0}^{t}\int_{\mathbb{R}_{*}^{2}}z\tilde{\Pi}\left(dt,dz\right)$, for $t\in \left[ 0,T \right]$,
is given by 
\[
\mathbb{V}\left[L\left(t\right)\right]=t\cdot\delta\frac{\Gamma\left(2-\alpha\right)}{\lambda^{2-\alpha}}\left(\rho\lambda^{1-\alpha}+\eta\eta^{\top}\right).
\]
We refer to Appendix \ref{app: Normal Tempered Stable process} for
further details.

Table \ref{tab: test - Parameters sets used in the numerical experiment}
\begin{table}
\caption{\label{tab: test - Parameters sets used in the numerical experiment}Parameter
sets}
\smallskip{}
\begin{centering}
\begin{tabular}{ccccc}
\toprule 
Parameters & VG0 & VG1 & NIG0 & NIG1\tabularnewline
\midrule
\midrule 
$\alpha$ & 0 & 0 & $\frac{1}{2}$ & $\frac{1}{2}$\tabularnewline
\midrule 
$\lambda$ & 1 & 6 & 20766.4 & 57.1108\tabularnewline
\midrule 
$\delta$ & 1 & 6 & 0.77576 & 4.26367\tabularnewline
\midrule 
$\eta^{\left(1\right)}$ & -0.1 & -0.1 & -37.688 & -0.295846\tabularnewline
\midrule 
$\eta^{\left(2\right)}$ & -0.2 & -0.2 & -2.224 & -0.292984\tabularnewline
\midrule 
$\rho^{\left(1,1\right)}$ & 0.09 & 0.01 & 3.984 & 0.037021\tabularnewline
\midrule 
$\rho^{\left(1,2\right)}$ & 0.06 & 0 & 3.160 & 0.026574\tabularnewline
\midrule 
$\rho^{\left(2,2\right)}$ & 0.16 & 0.0225 & 3.512 & 0.054613\tabularnewline
\midrule 
$r$ & 0.05 & 0 & 0 & 0\tabularnewline
\midrule 
$T$ & 1 & $\frac{1}{2}$ & $\frac{1}{2}$ & $\frac{1}{2}$\tabularnewline
\midrule 
$K$ & 100 & 100 & 100 & 100\tabularnewline
\bottomrule
\end{tabular}
\par\end{centering}
\end{table}
 lists four sets of representative parameter values where we always take the diffusion matrix $\sigma$ equal to zero. Table \ref{tab: test - Variance of the parameters sets}
\begin{table}
\caption{\label{tab: test - Variance of the parameters sets}Standard deviation and correlation coefficient}
\smallskip{}
\begin{centering}
\begin{tabular}{ccccc}
\toprule 
 & VG0 & VG1 & NIG0 & NIG1\tabularnewline
\midrule
\midrule 
$\sqrt{\mathbb{V}\left[L^{\left(1\right)}\left(1\right)\right]}$ & 0.3162 & 0.1080 & 0.1958 & 0.1943\tabularnewline
\midrule 
$\sqrt{\mathbb{V}\left[L^{\left(2\right)}\left(1\right)\right]}$ & 0.4472 & 0.1707 & 0.1830 & 0.2352\tabularnewline
\midrule 
$\frac{\textrm{cov}\left[L^{\left(1\right)}\left(1\right),L^{\left(2\right)}\left(1\right)\right]}{\sqrt{\mathbb{V}\left[L^{\left(1\right)}\left(1\right)\right]\mathbb{V}\left[L^{\left(2\right)}\left(1\right)\right]}}$ & 0.5656 & 0.1807 & 0.8417 & 0.5975\tabularnewline
\bottomrule
\end{tabular}\medskip{}
\par\end{centering}
\end{table}
 contains the corresponding standard deviations and correlation coefficients.
The sets VG0 and NIG0 are taken from \citet[page 208]{hilber2013}
and \citet[Figure 8]{rydberg1997}, respectively. The VG1 set was
designed by us based on VG0. Finally, the NIG1 set was obtained via
standard maximum likelihood estimation\footnote{The density function for the case where $\alpha\in\left\{ 0,\frac{1}{2}\right\} $
can be found in Appendix \ref{app: Normal Tempered Stable process}.} using the close price data of S\&P500 (\textasciicircum GSPC) and
EUROSTOXX50 (\textasciicircum STOXX50E), retrieved from Yahoo Finance,
covering the period from 01/01/2014 to 31/12/2024. In particular,
we apply the methodology used by \citet{hainaut2014} to the logarithmic
return of the price, i.e. $d\ln X$.

The following list specifies all choices for the values of the parameters of our numerical scheme:
\begin{itemize}
\item $N_{z}=2N_{x}$ and $N_{t}=\text{round}\left[\frac{1}{2}N_{x}\right]$.
Clearly, with this choice, the three mesh widths are directly proportional
to each other.
\item $z_{\max}^{\mathbf{I}}=2h_{z}$. This choice is motivated by the fact
that the artificial diffusion acts over a small region around the
origin. 
\item $z_{\max}^{\mathbf{II}}=h_z \cdot \text{ceil}\left[\frac{\sqrt{0.1}z_{\max}^{\mathbf{III}}}{h_z}-\frac{1}{2}\right]$, such that\footnote{Whereas $z_{\max}^{\mathbf{II}}$ is not independent of $h_z$, the dependence is weak and the proof in Appendix \ref{sec: appendix - Modified midpoint rule} is readily adapted such that the convergence result of Proposition~\ref{prop: integral disc - Quadrature scheme for the integral operator} remains valid.} $z_{\max}^{\mathbf{II}}\rightarrow\sqrt{0.1}z_{\max}^{\mathbf{III}}$ as $h_z \rightarrow 0^{+}$. This means that the size of $R_{z}^{\mathbf{II}}$ is about 10\% of the full integration domain $R_{z}$.
\item $z_{\max}^{\mathbf{III}}=\max \left\{ \left\Vert z\right\Vert _{\infty}: z\in\mathbb{R}^{2}, \ell\left(z\right)=10^{-8}\right\} $.
Since the L{\'e}vy measure decays at least exponentially as $\left\Vert z \right\Vert \rightarrow\infty$,
we ensure that $\ell\left(z\right)<10^{-8}$ for all $z\in\mathbb{R}^{2}$
such that $\left\Vert z\right\Vert _{\infty}>z_{\max}^{\mathbf{III}}$. Thus, we choose $z_{\max}^{\mathbf{III}}$ equal to 11.5010 for VG0, 2.1410 for VG1, 0.4172 for NIG0, and 0.8807 for NIG1.
\item $x_{\rm int}=\frac{5}{2}K$. The non-smoothness set of $\phi$
is contained in the portion of $R_{x}$ where the grid $\mathbf{x}$ is uniform.
\item $x_{\rm max}$ was heuristically chosen equal to $57K$ for VG0, $5K$ for VG1, $6K$ for NIG0, and $7K$ for NIG1.
\item $c$ is chosen such that $F=\max\left(65\%,\frac{x_{\rm int}}{x_{\rm max}}\right)$.
In this way, approximately at least 65\% of the spatial grid points in each given direction belong to the interval $\left[0,x_{\rm int}\right]$. Thus, we choose $c$ equal to 21.6164 for VG0, 67.1487 for VG1, 55.0673 for NIG0, and 45.0189 for NIG1.
\item $N_{y}^{-}=\text{ceil}\left[-\frac{1}{h_{z}}\ln\left(x_{1}\right)\right]+N_{y}^{*}$
and $N_{y}^{+}=\text{ceil}\left[\frac{1}{h_{z}}\ln\left(x_{\rm max}\right)\right]+N_{y}^{*}$
for some given $N_{y}^{*}\in\mathbb{N}$. This choice minimizes the need for extrapolation in 
(\ref{eq: spatial disc - interpolation output value})
as it is necessary to extrapolate just to the grid points $x_{m_{1}m_{2}}$ with either $m_1=0$ or $m_2=0$.
This is done in a linear fashion.
In (\ref{eq: spatial disc - interpolation input value}),  as well as (\ref{eq: spatial disc - summation operator (2) over the y-grids}), we set $v\left(x,t\right)=0$ whenever $x\notin R_x$.
\item $N_{y}^{*}$ is taken as the minimal $n\in\mathbb{N}$ such that
the maximal prime factor of $\sharp \textrm{in} = N_{y}^{-}+N_{y}^{+}+2N_{z}$ is at most 7. This is beneficial for the efficiency of the FFT.
\item The tolerances for the fixed-point iteration and BiCGSTAB are
set to $10^{-7}$ and $10^{-14}$, respectively.
\item Interpolation is performed by cubic Lagrange polynomials. 
\end{itemize}

\begin{figure}
\caption{\label{fig: test - Surface solution}European put-on-the-average option
price and the Greeks Delta and Gamma for the parameter set NIG0}
\bigskip
\noindent \begin{centering}
\includegraphics[width=7.5cm]{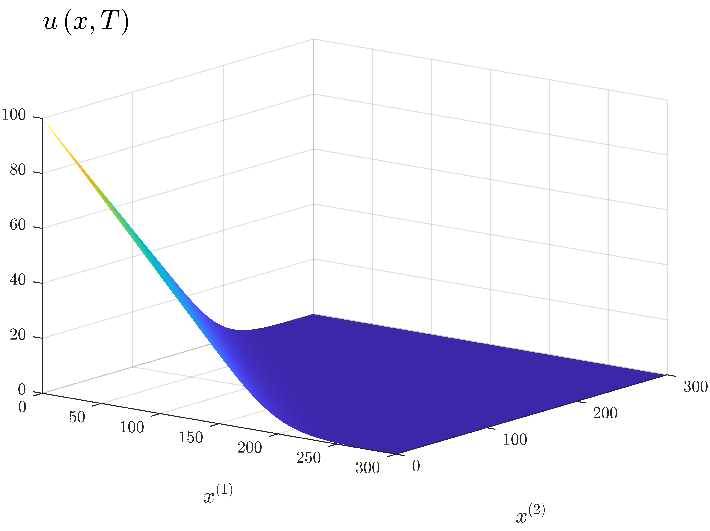}\hspace{0.5cm}\includegraphics[width=7.5cm]{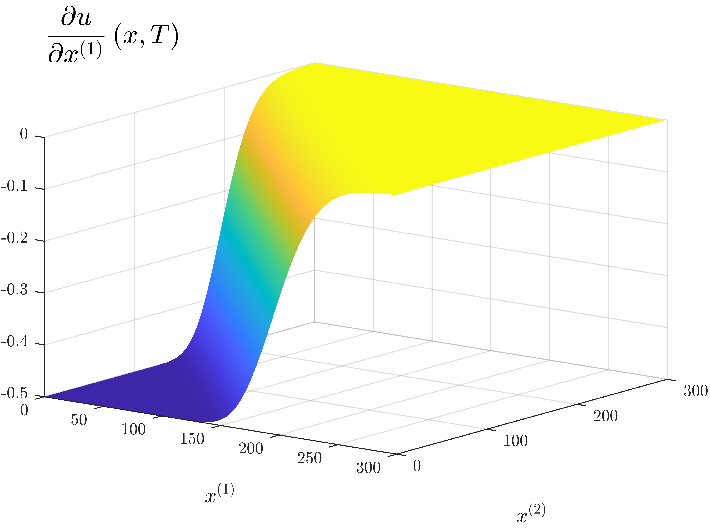}
\par\end{centering}
\noindent \begin{centering}
\includegraphics[width=7.5cm]{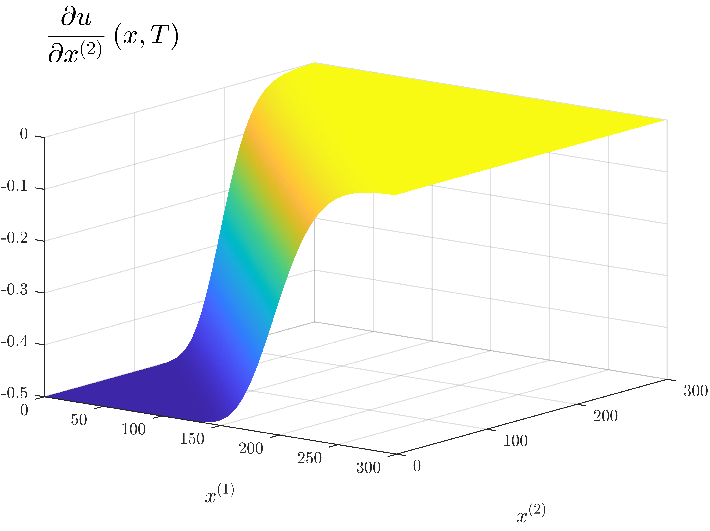}\hspace{0.5cm}\includegraphics[width=7.5cm]{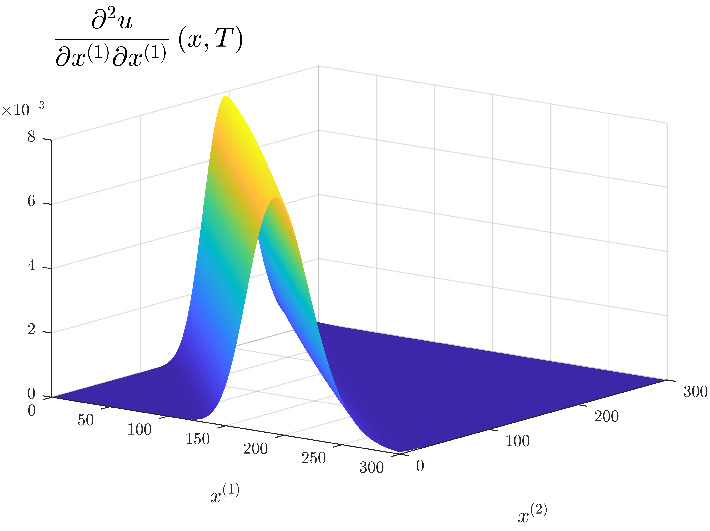}
\par\end{centering}
\noindent \begin{centering}
\includegraphics[width=7.5cm]{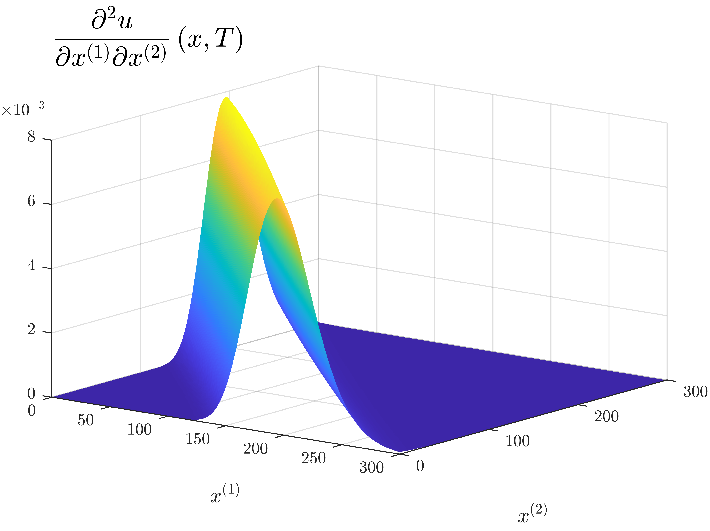}\hspace{0.5cm}\includegraphics[width=7.5cm]{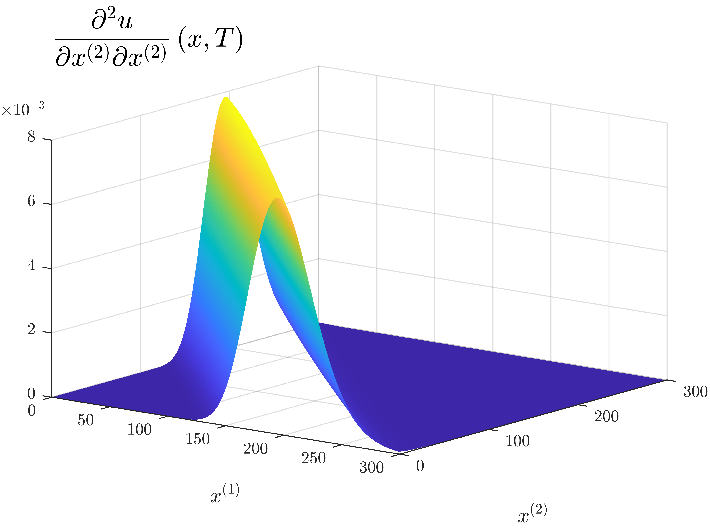}
\par\end{centering}
\end{figure}

Figure \ref{fig: test - Surface solution}
 displays the graphs of the option price function and its Greeks Delta and Gamma for the parameter set NIG0
from Table \ref{tab: test - Parameters sets used in the numerical experiment}, where we have taken $N_{x}=400$. The Greeks have been approximated (at negligible computational cost) by applying the second-order central finite difference schemes described in Section \ref{subsec: spatial discretization}.
Table \ref{tab: test - values of the reference solution}
\begin{table}
\caption{\label{tab: test - values of the reference solution}Numerical option prices for points $x$ near $\left(K,K\right)$}

\centering{}%
\begin{tabular}{ccccc}
\toprule 
$\left(x^{\left(1\right)},x^{\left(2\right)}\right)$ & VG0 & VG1 & NIG0 & NIG1\tabularnewline
\midrule
\midrule 
$\left(90,90\right)$ & 12.6540 & 10.1079 & 11.4067 & 11.5833\tabularnewline
\midrule 
$\left(90,100\right)$ & 10.6127 & 5.8453 & 7.8724 & 8.1532\tabularnewline
\midrule 
$\left(90,110\right)$ & 9.0142 & 3.0171 & 5.1023 & 5.4661\tabularnewline
\midrule 
$\left(100,90\right)$ & 10.4066 & 5.7627 & 7.8897 & 8.0913\tabularnewline
\midrule 
$\left(100,100\right)$ & 8.8020 & 2.9029 & 5.1186 & 5.3956\tabularnewline
\midrule 
$\left(100,110\right)$ & 7.5314 & 1.3889 & 3.1156 & 3.4314\tabularnewline
\midrule 
$\left(110,90\right)$ & 8.6186 & 2.8062 & 5.1393 & 5.3384\tabularnewline
\midrule 
$\left(110,100\right)$ & 7.3468 & 1.3058 & 3.1326 & 3.3739\tabularnewline
\midrule 
$\left(110,110\right)$ & 6.3294 & 0.6012 & 1.7937 & 2.0401\tabularnewline
\bottomrule
\end{tabular}
\end{table}
provides the numerical option prices for various points $x$ around $\left(K,K\right)$ and all parameter sets from Table \ref{tab: test - Parameters sets used in the numerical experiment}. Here $N_x=800$ for VG0, VG1, and NIG1, while for NIG0 we used $N_x=400$.

\begin{figure}
\caption{\label{fig: test - Error in ROI}Total error in $\left[0,3K\right]\times\left[0,3K\right]$}

\noindent \centering{}\includegraphics[width=16cm]{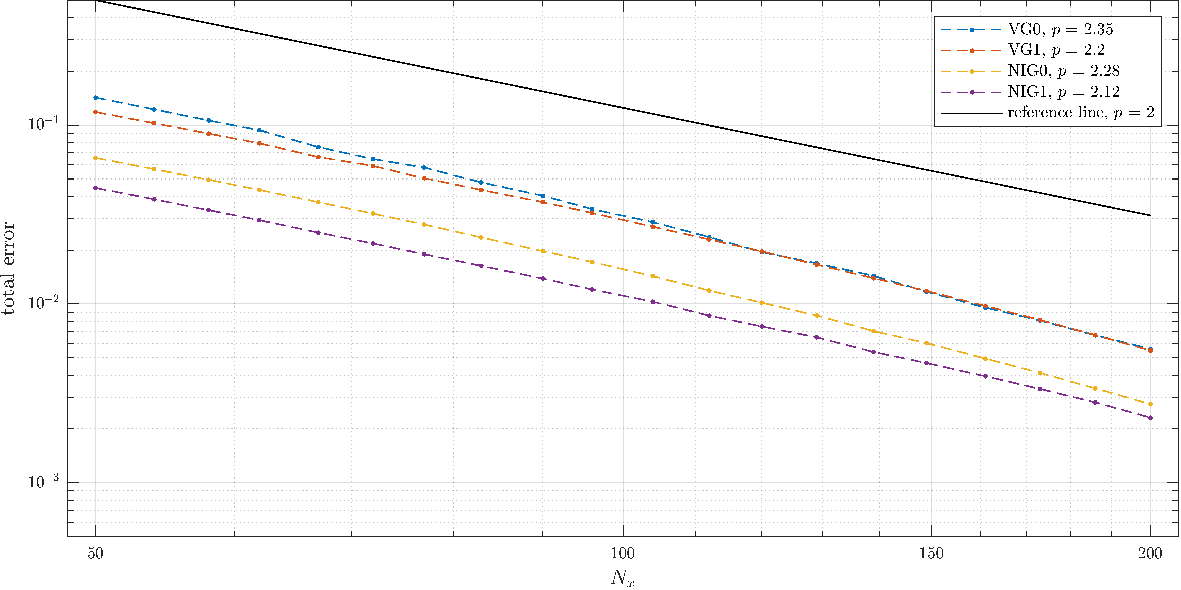}
\end{figure}

We next investigate the convergence behaviour of the numerical scheme. 
Let $\mathbf{x}_N$ denote the set of spatial grid points if $N_{x} = N$. For $x\in R_{x}$, let $\tilde{u}(x;N)$ denote the approximation of the exact solution value $u(x,T)$ obtained by the numerical scheme if $N_{x} = N$. More precisely, the vector $V^{N_t}$ generated by \eqref{eq: temporal disc - fixed-point (1) recursive equation}-\eqref{eq: temporal disc - fixed-point (3) exit condition} yields the approximation on the spatial grid $\mathbf{x}_N$ and cubic interpolation is employed whenever $x\notin \mathbf{x}_N$.
We consider $\tilde{u}(x;N)$ with $N=400$ as the reference solution and study for $50\leq N \leq 200$ the total error defined by
\[
E\left(N\right)=\max\left\{ |\tilde{u}(x;N)-\tilde{u}(x;400)|: x\in \mathbf{x}_N \textrm{~and~} x\in\left[0,3K\right]\times\left[0,3K\right]\right\} .
\]
Figure \ref{fig: test - Error in ROI}
displays the total errors for all four parameter sets from Table \ref{tab: test - Parameters sets used in the numerical experiment}.
The quantity $p$ in the legend denotes the numerical order of convergence, which is computed by linear regression.
Clearly, the favourable result is found that the numerical scheme achieves second-order convergence for each set of parameters. 
This is in agreement with the theoretical order of convergence for the discretization of the integral term given by Proposition \ref{prop: integral disc - Quadrature scheme for the integral operator} in the case of the two VG sets (where $A_{\ell}=0$).
On the other hand, in the case of the two NIG sets (where $A_{\ell}=1$), the theoretical order of convergence is only one. 
At present we have no clear explanation why, for the two NIG sets, second-order convergence is observed in the numerical experiments. We shall leave this interesting question for future research.

\section{\label{sec: Conclusions}Conclusions}

In this paper, we have developed an effective numerical method for the valuation of European options under two-asset exponential L{\'e}vy models with particular attention to the infinite-activity case. 

Our method is based upon the ideas in \citet{wang2007} for the one-asset case.
A key part of our method is the tailored discretization of the non-local double integral term, designed to handle singular measures under mild assumptions. The discretized integral term can subsequently be efficiently evaluated by FFT.
For the discretization in time, the semi-Lagrangian Crank--Nicolson method is employed with a fixed-point iteration on the integral part.

Numerical experiments for put-on-the-average options under Normal Tempered Stable processes indicate that our method achieves favourable second-order convergence whenever the exponential L{\'e}vy model has finite-variation. A relevant rigorous theoretical result is proved on the convergence behaviour of the discretization of the integral term.

The numerical method derived in this paper is easily adapted to many other European options on two assets, such as spread options in energy markets.

A main topic for future research will be extending the proposed methodology to the valuation of American-style two-asset options under exponential L{\'e}vy models with infinite-activity, where the combination of the early-exercise feature and the non-local integral term poses additional challenges.

\section{Acknowledgements}

The authors acknowledge the support of the Research Foundation - Flanders (FWO) under grant G0B5623N and the FWO Scientific Research Network
ModSimFIE (FWO WOG W001021N).  The third author also acknowledges the financial support of the Research Foundation - Flanders (FWO) through FWO SAB K803124.

\newpage{}

\appendix

\section{\label{app: Normal Tempered Stable process}$d$-dimensional Normal
Tempered Stable process}

The term $d$-dimensional Normal Tempered Stable process refers to a
$d$-dimensional pure-jump compensated L{\'e}vy process $L$ with L{\'e}vy
measure generated by subordinating a $d$-dimensional Brownian motion
$B$ with a tempered stable subordinator $G$, i.e., a pure-jump process
with almost surely non-decreasing trajectories. Such a process is defined
by the following equation
\[
L\left(t\right)=B\left(G\left(t\right)\right)-\mathbb{E}\left[B\left(G\left(t\right)\right)\right]\qquad\text{with }L\left(0\right)=0.
\]
In our context, we will use this process to define the jump component
of the logarithmic return in asset prices, i.e., we choose 
\[
\int_{\mathbb{R}_{*}^{2}}z\tilde{\Pi}\left(dt,dz\right)=dL\left(t\right).
\]

\subsection{Tempered Stable subordinator}

A tempered stable subordinator is a non-compensated $1$-sided tempered
stable process $G$, which is characterized by the parameters  $\delta, \lambda>0$ and $\alpha\in\left[0,1\right)$. For more details see \citet{kuchler2013}. Table~\ref{tab: TS subordinator}
\begin{table}
\caption{\label{tab: TS subordinator}Main quantities of the Tempered Stable subordinator}

\smallskip{}

\noindent \begin{centering}
\begin{tabular}{cc}
\toprule 
Quantity & Formula\tabularnewline
\midrule
\midrule 
L{\'e}vy measure & $\ell_{G}\left(z\right)=\mathbb{I}_{z>0}\delta e^{-\lambda z}z^{-1-\alpha}$\tabularnewline
\midrule 
Characteristic exponent & $\psi_{G}\left(z\right)=\begin{cases}
-\delta\ln\left(1-iz\lambda^{-1}\right) & \text{if }\alpha=0\\
\delta\Gamma\left(-\alpha\right)\left(\left(\lambda-iz\right)^{\alpha}-\lambda^{\alpha}\right) & \text{if }\alpha\in\left(0,1\right)
\end{cases}$\tabularnewline
\midrule 
Expected value & $\mathbb{E}\left[G\left(1\right)\right]=\delta\frac{\Gamma\left(1-\alpha\right)}{\lambda^{1-\alpha}}$\tabularnewline
\midrule 
Variance & $\mathbb{V}\left[G\left(1\right)\right]=\delta\frac{\Gamma\left(2-\alpha\right)}{\lambda^{2-\alpha}}$\tabularnewline
\midrule 
Density function & $f_{G}\left(z\right)=\begin{cases}
\mathbb{I}_{z>0}\frac{\lambda^{\delta}}{\Gamma\left(\delta\right)}z^{\delta-1}e^{-\lambda z} & \text{if }\alpha=0\\
\mathbb{I}_{z>0}\delta z^{-\frac{3}{2}}e^{-\left(\sqrt{\lambda}z-\delta\sqrt{\pi}\right)^{2}z^{-1}} & \text{if }\alpha=\frac{1}{2}\\
\text{not known analytically} & \text{else}
\end{cases}$\tabularnewline
\bottomrule
\end{tabular}\smallskip{}
\par\end{centering}
\end{table}
 shows the main quantities for such a process. Note that $G$ corresponds to the Gamma process for $\alpha=0$ and to the Inverse Gaussian process for $\alpha=\frac{1}{2}$.

\subsection{Normal Tempered Stable process}

Consider $B\left(t\right)=\eta t+\sqrt{\rho}\,W\left(t\right)$, where $W$  is a standard $d$-dimensional Wiener process,  $\eta\in\mathbb{R}^{d}$  and  $\sqrt{\rho}$  is obtained by the Cholesky decomposition of a given positive definite symmetric matrix $\rho$, i.e. $\rho=\sqrt{\rho}\cdot\sqrt{\rho}^{\top}$. Adapting the results presented in \citet{barndorff-nielsen2001} and \citet[Chapter 4]{rocha-arteaga2019}\footnote{The authors consider the more general case where the characteristic
exponent of $L$ is defined as $\psi_{L}\left(\tau\right)=\int_{\mathbb{R}^{d}}\left(e^{i\tau^{\top}z}-1-i\tau^{\top}z\mathbb{I}_{\left\Vert z\right\Vert <1}\right)\ell_{L}\left(dz\right)$,
while we consider the case where $\psi_{L}\left(\tau\right)=\int_{\mathbb{R}^{d}}\left(e^{i\tau^{\top}z}-1-i\tau^{\top}z\right)\ell_{L}\left(dz\right)$.}, we define a $d$-dimensional Normal Tempered Stable process as 
\begin{equation}
L\left(t\right)=B\left(G\left(t\right)\right)-ct,\label{eq: NTS: definition}
\end{equation}
where $G$ is a Tempered Stable subordinator and $c=\mathbb{E}\left[B\left(G\left(t\right)\right)\right]=\delta\frac{\Gamma\left(1-\alpha\right)}{\lambda^{1-\alpha}}\eta$.

Table \ref{tab: NTS}
\begin{table}
\caption{\label{tab: NTS}Main quantities of the Normal Tempered Stable process}

\smallskip{}

\noindent \begin{centering}
\begin{tabular}{cc}
\toprule 
Quantity & Formula\tabularnewline
\midrule
\midrule 
L{\'e}vy measure & $\ell_{L}\left(z\right)=\delta\Phi\left(z\mid\alpha,0\right)$\tabularnewline
\midrule 
Characteristic exponent & $\psi_{L}\left(z\right)=\begin{cases}
-\delta\ln\left(\frac{\lambda-iz^{\top}\eta+\frac{1}{2}z^{\top}\rho z}{\lambda}\right)-iz^{\top}c & \text{if }\alpha=0\\
\delta\Gamma\left(-\alpha\right)\left(\left(\lambda-iz^{\top}\eta+\frac{1}{2}z^{\top}\rho z\right)^{\alpha}-\lambda^{\alpha}\right)-iz^{\top}c & \text{if }\alpha\in\left(0,1\right)
\end{cases}$\tabularnewline
\midrule 
Expected value & $\mathbb{E}\left[L\left(1\right)\right]=0$\tabularnewline
\midrule 
Variance & $\mathbb{V}\left[L\left(1\right)\right]=\delta\frac{\Gamma\left(2-\alpha\right)}{\lambda^{2-\alpha}}\left(\rho\frac{\lambda}{1-\alpha}+\eta\eta^{\top}\right)$\tabularnewline
\midrule 
Density function & $f_{L}\left(z\right)=\begin{cases}
\frac{\lambda^{\delta}}{\Gamma\left(\delta\right)}\Phi\left(z+c\mid-\delta,0\right) & \text{if }\alpha=0\\
\delta e^{2\delta\sqrt{\lambda\pi}}\Phi\left(z+c\mid\frac{1}{2},\delta^{2}\pi\right) & \text{if }\alpha=\frac{1}{2}\\
\text{not known analytically} & \text{else}
\end{cases}$\tabularnewline
\bottomrule
\end{tabular}\smallskip{}
\par\end{centering}
\end{table}
 shows the main quantities for such a process. Most of the formulae
are expressed in terms of the function $\Phi$ which is given by
\begin{equation}
\Phi\left(z\mid a,b\right)=2\sqrt{\frac{\left(\left\Vert \eta\right\Vert _{\rho}^{2}+2\lambda\right)^{a+\frac{d}{2}}}{\left(2\pi\right)^{d}\det\left[\rho\right]}}\frac{K_{a+\frac{d}{2}}\left(\sqrt{\left(\left\Vert \eta\right\Vert _{\rho}^{2}+2\lambda\right)\left(\left\Vert z\right\Vert _{\rho}^{2}+2b\right)}\right)}{\left(\sqrt{\left\Vert z\right\Vert _{\rho}^{2}+2b}\right)^{a+\frac{d}{2}}}e^{\left\langle \eta,z\right\rangle _{\rho}}\label{eq: NTS, function Phi}
\end{equation}
where $K_{\nu}\left(\tau \right)=\frac{1}{2}\int_{0}^{\infty}y^{\nu-1}e^{-\frac{1}{2}\tau  \left(y+y^{-1}\right)}dy$, for $\tau>0$,
denotes the modified Bessel function of the second kind (see \citet[Appendix A]{schoutens2003}),
$\left\langle x,y\right\rangle _{\rho}=x^{\top}\rho^{-1}y$ and $\left\Vert x\right\Vert _{\rho}=\sqrt{\left\langle x,x\right\rangle _\rho}$ is its induced norm.  We conclude this appendix with the following proposition. 

\begin{prop}
Consider a L{\'e}vy measure $\ell$ over $\mathbb{R}_{*}^{d}=\mathbb{R}^{d}\setminus\left\{ 0\right\} $.
Assume that there exist constants $A_{\ell}$ and $B_{\ell}$, and
for any given $h>0$ two constants $C_{\ell}\left(h\right)$ and $C_{\ell}^{\prime}\left(h\right)$
such that

\[
\begin{cases}
\ell\left(z\right)\leq C_{\ell}\left(h\right)\left\Vert z\right\Vert_{\rho}^{-2-A_{\ell}} & \text{for any }z\text{ such that }\left\Vert z\right\Vert_{\rho}\in\left(0,h\right],\\

\left|\ell_{z}^{\left(j\right)}\left(z\right)\right|\leq C_{\ell}^{\prime}\left(h\right)\left\Vert z\right\Vert_{\rho} ^{-3-A_{\ell}}& \text{for any }z\text{ such that }\left\Vert z\right\Vert_{\rho}\in\left(0,h\right],\\

\ell\left(z\right)=O\left(e^{-B_{\ell}\left\Vert z\right\Vert_{\rho} }\right) & \text{as }\left\Vert z\right\Vert_{\rho} \rightarrow\infty,
\end{cases}
\]
where $\ell_{z}^{\left(j\right)}$ denotes the $j$-th partial derivative of $\ell$ with respect to $z$, with $j=1,2,\ldots,d$.\\
Then, for a Normal Tempered Stable process these constants are given by 
\begin{align*}
& A_{\ell}  =2\alpha,\qquad 
B_{\ell}  =\sqrt{\left\Vert \eta\right\Vert _{\rho}^{2}+2\lambda}-\left\Vert \eta\right\Vert _{\rho},\\
&C_{\ell}\left(h\right) =\delta\frac{2^{\alpha+\frac{d}{2}}\Gamma\left(\alpha+\frac{d}{2}\right)}{\sqrt{\left(2\pi\right)^{d}\det\left[\rho\right]}}e^{h\left\Vert \eta\right\Vert _{\rho}}\, ,\\
& C_{\ell}^{\prime}\left(h\right)  =C_{\ell}\left(h\right)\left\Vert \rho^{\frac{1}{2}}\right\Vert _{\rho}\left\Vert \rho^{-1}\right\Vert _{\rho}\left(h\left\Vert \eta\right\Vert _{\rho}+\left(A_{\ell}+d\right)\right).
\end{align*}
\end{prop}

\begin{proof}
Defining $c_{1}=\sqrt{\left\Vert \eta\right\Vert _{\rho}^{2}+2\lambda}$,
$c_{2}=2\delta c_{1}^{2\alpha+d}\left(2\pi\right)^{-\frac{d}{2}}\det\left[\rho\right]^{-\frac{1}{2}}$
and $\hat{K}_{\nu}\left(\tau\right)=\tau^{-\nu}K_{\nu}\left(\tau\right)$,
the L{\'e}vy measure and its gradient can be expressed as 
\begin{align*}
\ell\left(z\right) & =c_{2}\hat{K}_{\alpha+\frac{d}{2}}\left(c_{1}\left\Vert z\right\Vert _{\rho}\right)e^{\left\langle \eta,z\right\rangle _{\rho}},\\
\ell_{z}\left(z\right) & =c_{2}e^{\left\langle \eta,z\right\rangle _{\rho}}\rho^{-1}\left(\hat{K}_{\alpha+\frac{d}{2}}\left(c_{1}\left\Vert z\right\Vert _{\rho}\right)\eta-\hat{K}_{\alpha+\frac{d}{2}+1}\left(c_{1}\left\Vert z\right\Vert _{\rho}\right)c_{1}^{2}z\right),
\end{align*}
where we have used the well known formulae $\frac{\partial}{\partial z}\left\langle \eta,z\right\rangle _{\rho}=\rho^{-1}\eta$,
$\frac{\partial}{\partial z}\left\Vert z\right\Vert _{\rho}=\left\Vert z\right\Vert _{\rho}^{-1}\rho^{-1}z$
and $\hat{K}_{\nu}^{\prime}\left(\tau\right)=-\tau\hat{K}_{\nu+1}\left(\tau\right)$.
By the Cauchy--Schwarz inequality, i.e. $\left\langle \eta,z\right\rangle _{\rho}\leq\left\Vert \eta\right\Vert _{\rho}\left\Vert z\right\Vert _{\rho}$,
and the asymptotic behaviour of $K_{\nu}$ as $\tau\rightarrow+\infty$,
i.e. $K_{\nu}\left(\tau\right)=O\left(\tau^{-\frac{1}{2}}e^{-\tau}\right)$,
we get that 
\[
\ell\left(z\right)=O\left(\left\Vert z\right\Vert _{\rho}^{-\left(\alpha+\frac{d+1}{2}\right)}e^{-\left(c_{1}-\left\Vert \eta\right\Vert _{\rho}\right)\left\Vert z\right\Vert _{\rho}}\right).
\]
Since $\left\Vert z\right\Vert _{\rho}^{-\left(\alpha+\frac{d+1}{2}\right)}=O\left(1\right)$
as $\left\Vert z\right\Vert _{\rho}\rightarrow+\infty$ for any $\alpha\in\left[0,1\right)$,
then we deduce that 
\[
B_{\ell}=\sqrt{\left\Vert \eta\right\Vert _{\rho}^{2}+2\lambda}-\left\Vert \eta\right\Vert _{\rho},
\]
with $B_{\ell}$ being positive for any choice of parameters. Since
from a well known bound of $K_{\nu}$ follows 
\[
\hat{K}_{\nu}\left(\tau\right)\leq2^{\nu-1}\Gamma\left(\nu\right)\tau^{-2\nu}\qquad\text{for any }\tau,\nu>0,
\]
then the L{\'e}vy measure is bounded by 
\[
\ell\left(z\right)\leq\frac{c_{2}}{2}\left(\frac{\sqrt{2}}{c_{1}}\right)^{2\alpha+d}\Gamma\left(\alpha+\frac{d}{2}\right)e^{\left\Vert \eta\right\Vert _{\rho}\left\Vert z\right\Vert _{\rho}}\left\Vert z\right\Vert _{\rho}^{-2\alpha-d}.
\]
Hence, 
\[
A_{\ell}=2\alpha,\qquad C_{\ell}\left(h\right)=\delta\frac{2^{\alpha+\frac{d}{2}}\Gamma\left(\alpha+\frac{d}{2}\right)}{\sqrt{\left(2\pi\right)^{d}\det\left[\rho\right]}}e^{h\left\Vert \eta\right\Vert _{\rho}}.
\]
Regarding the condition on the gradient of $\ell$, we get similarly
the following bound
\begin{align*}
\left\Vert \ell_{z}\left(z\right)\right\Vert _{\rho} & \leq c_{2}e^{h\left\Vert \eta\right\Vert _{\rho}}\left\Vert \rho^{-1}\right\Vert _{\rho}\left(\hat{K}_{\alpha+\frac{d}{2}}\left(c_{1}\left\Vert z\right\Vert _{\rho}\right)\left\Vert \eta\right\Vert _{\rho}+\hat{K}_{\alpha+\frac{d}{2}+1}\left(c_{1}\left\Vert z\right\Vert _{\rho}\right)c_{1}^{2}\left\Vert z\right\Vert _{\rho}\right)\\
 & \leq C_{\ell}\left(h\right)\left\Vert \rho^{-1}\right\Vert _{\rho}\left(\left\Vert \eta\right\Vert _{\rho}\left\Vert z\right\Vert _{\rho}^{-A_{\ell}-d}+\left(A_{\ell}+d\right)\left\Vert z\right\Vert _{\rho}^{-A_{\ell}-d-1}\right)\\
 & \leq C_{\ell}\left(h\right)\left\Vert \rho^{-1}\right\Vert _{\rho}\left(h\left\Vert \eta\right\Vert _{\rho}+\left(A_{\ell}+d\right)\right)\left\Vert z\right\Vert _{\rho}^{-A_{\ell}-d-1}.
\end{align*}
Consider now that $\left|\ell_{z}^{\left(i\right)}\left(z\right)\right|\leq\left\Vert \ell_{z}\left(z\right)\right\Vert _{2}$
and that for any vector $z\in\mathbb{R}^{d}$
\[
\left\Vert z\right\Vert _{2}=z^{\top}z=z^{\top}\rho^{\frac{1}{2}}\rho^{-1}\rho^{\frac{1}{2}}z=\left(\rho^{\frac{1}{2}}z\right)^{\top}\rho^{-1}\left(\rho^{\frac{1}{2}}z\right)=\left\Vert \rho^{\frac{1}{2}}z\right\Vert _{\rho}\leq\left\Vert \rho^{\frac{1}{2}}\right\Vert _{\rho}\left\Vert z\right\Vert _{\rho}.
\]
Then, 
\[
\left|\ell_{z}^{\left(i\right)}\left(z\right)\right|\leq C_{\ell}\left(h\right)\left\Vert \rho^{\frac{1}{2}}\right\Vert _{\rho}\left\Vert \rho^{-1}\right\Vert _{\rho}\left(h\left\Vert \eta\right\Vert _{\rho}+\left(A_{\ell}+d\right)\right)\left\Vert z\right\Vert _{\rho}^{-A_{\ell}-d-1},
\]
from which we deduce that 
\[
C_{\ell}^{\prime}\left(h\right)=C_{\ell}\left(h\right)\left\Vert \rho^{\frac{1}{2}}\right\Vert _{\rho}\left\Vert \rho^{-1}\right\Vert _{\rho}\left(h\left\Vert \eta\right\Vert _{\rho}+\left(A_{\ell}+d\right)\right).
\]
\end{proof}

\section{\label{sec: appendix - Modified midpoint rule}Proof of Proposition \ref{prop: integral disc - Quadrature scheme for the integral operator}}
Write $S=R_x \times\left[0,T\right]$. To analyze for $\left(x,t\right)\in S$ the order of convergence of $E\left(x,t\right)$ in (\ref{eq: error formula}), we decompose the error according to the three approximations:
\begin{equation}
E\left(x,t\right)=E^{\left(1\right)}\left(x,t\right)+E^{\left(2\right)}\left(x,t\right)+E^{\left(3\right)}\left(x,t\right),\label{eq: errore decomposition}
\end{equation}
where 
\begin{align*}
E^{\left(1\right)}\left(x,t\right) & =\int_{R_{z}^{\mathbf{I}}}f\left(z,x,t\right)\ell\left(dz\right)-\frac{1}{2}\mathbf{1}^{\top}\left(u_{xx}\left(x,t\right)\circ I_{x}\left(\int_{R_{z}^{\mathbf{I}}}zz^{\top}\ell\left(dz\right)\right)I_{x}\right)\mathbf{1},\\
E^{\left(2\right)}\left(x,t\right) & =\int_{R_{z}^{\mathbf{II}}}f\left(z,x,t\right)\ell\left(dz\right)-\sum_{l_{1},l_{2}:z_{l_{1}l_{2}}\in R_{z}^{\mathbf{II}}}\omega_{l_{1}l_{2}}f\left(z_{l_{1}l_{2}},x,t\right),\\
E^{\left(3\right)}\left(x,t\right) & =\int_{R_{z}^{\mathbf{III}}}f\left(z,x,t\right)\ell\left(dz\right)-\sum_{l_{1},l_{2}:z_{l_{1}l_{2}}\in R_{z}^{\mathbf{III}}}\omega_{l_{1}l_{2}}f\left(z_{l_{1}l_{2}},x,t\right).
\end{align*}
In what follows, for notational convenience, we shall omit the dependence on $\left(x,t\right)$. 

For the proof, we assume:
$u\left(\cdot,t\right)\in\mathcal{C}^{3}\left(\mathbb{R}_{\geq0}^{2}\right)$ whenever $t\in\left[0,T\right]$; $\ell\in\mathcal{C}^{2}\left(\mathbb{R}_{*}^{2}\right)$;
$z_{\max}^{\mathbf{I}}=Mh_{z}$ for some fixed integer $M\ge 1$. 
For any given $a,b,c\in\mathbb{R}$ such that $0\leq a\leq b$, the following useful bounds hold\footnote{Use that $\left\Vert z\right\Vert _{\infty}\leq\left\Vert z\right\Vert \leq\sqrt{2}\left\Vert z\right\Vert _{\infty}$ for any $z\in\mathbb{R}^{2}$ together with a conversion of the integral to polar coordinates.}
\begin{equation} 
\textrm{if } c\not= -2: \quad
\int_{a\leq\left\Vert z\right\Vert _{\infty}\leq b}\left\Vert z\right\Vert ^{c}dz \leq 
\frac{2\pi}{c+2} \left( (\sqrt{2}\,b)^{c+2}-a^{c+2}\right),
\label{eq: proof - radial integral}
\end{equation}
\begin{equation} 
\textrm{if } a>0: \quad
\int_{a\leq\left\Vert z\right\Vert _{\infty}\leq b}\left\Vert z\right\Vert ^{-2}dz \leq 
2\pi \left( \ln(\sqrt{2}\,b)-\ln(a)\right).
\label{eq: proof - radial integral added}
\end{equation}

Define the functions $\hat{f}\left(z\right)=f\left(z\right)\left\Vert z\right\Vert ^{-2}$
and $\hat{\ell}\left(z\right)=\ell\left(z\right)\left\Vert z\right\Vert ^{2}$ ($z\in \mathbb{R}_{*}^{2}$).
It can be seen that $\hat{f}\in\mathcal{C}^{2}\left(\mathbb{R}_{*}^{2}\right)$ and there exist constants $C_{f}^{\prime}$,
$C_{f}^{\prime\prime}$ such that\footnote{The property (\ref{eq: proof - bounds on derivatives}) is weaker than in the 1-dimensional case, where the first and second derivatives of $\hat{f}$ are uniformly bounded near $z=0$, see \citet{wang2007}.}
\begin{equation}
    \left|\hat{f}_{z}^{\left(i\right)}\left(z\right)\right| \le \frac{C_{f}^{\prime}}{\left\Vert z\right\Vert} 
\quad\textrm{and}\quad
\left|\hat{f}_{zz}^{\left(i,j\right)}\left(z\right)\right| \le \frac{C_{f}^{\prime\prime}}{\left\Vert z\right\Vert^2} 
\label{eq: proof - bounds on derivatives}
\end{equation}
whenever $i,j \in \{1,2\}$, $0< \left\Vert z\right\Vert _{\infty}\leq z_{\max}^{\mathbf{II}}$ and $\left(x,t\right)\in S$.

Start by $E^{\left(1\right)}$. From Section \ref{subsec: integral disc}, this error is given by 
\[
E^{\left(1\right)}=\int_{R_{z}^{\mathbf{I}}}\varepsilon\left(z\right)\ell\left(dz\right),
\]
where $\varepsilon$ satisfies $\varepsilon\left(z\right)=O\left(\left\Vert z\right\Vert ^{3}\right)$ uniformly in $\left(x,t\right)\in S$.
By the upper bound of the L{\'e}vy measure (\ref{eq:  intro - upper bound of the L=0000E9vy measure}) and the monotonicity of the integral, it follows that $E^{\left(1\right)}=O\left(\int_{R_{z}^{\mathbf{I}}}\left\Vert z\right\Vert ^{1-A_{\ell}}dz\right)$.
Hence, using (\ref{eq: proof - radial integral}) and recalling $A_{\ell}<2$, we deduce that
\begin{equation}
E^{\left(1\right)}=O\left(h_{z}^{3-A_{\ell}}\right)\label{eq: E^1}
\end{equation}
uniformly in $\left(x,t\right)\in S$.

Now, move on to $E^{\left(2\right)}$. Since $\hat{f}\in\mathcal{C}^{2}\left(R_{z}^{\mathbf{II}}\right)$, its Taylor approximation for $z\in R_{l_{1}l_{2}}$ is given
by 
\[
\hat{f}\left(z\right)=\hat{f}\left(z_{l_{1}l_{2}}\right)+\left(z-z_{l_{1}l_{2}}\right)^{\top}\hat{f}_{z}\left(z_{l_{1}l_{2}}\right)+\frac{1}{2}\left(z-z_{l_{1}l_{2}}\right)^{\top}\hat{f}_{zz}\left(\zeta_{l_{1}l_{2}}\left(z\right)\right)\left(z-z_{l_{1}l_{2}}\right),
\]
where $\zeta_{l_{1}l_{2}}\left(z\right)\in R_{l_{1}l_{2}}$ depends
on $z$. By substitution, this yields
\begin{align}
E^{\left(2\right)} & =\sum_{l_{1},l_{2}:z_{l_{1}l_{2}}\in R_{z}^{\mathbf{II}}}\int_{R_{l_{1}l_{2}}}\left(\hat{f}\left(z\right)-\hat{f}\left(z_{l_{1}l_{2}}\right)\right)\hat{\ell}\left(z\right)dz\nonumber \\
 & =\sum_{l_{1},l_{2}:z_{l_{1}l_{2}}\in R_{z}^{\mathbf{II}}}\left(\int_{R_{l_{1}l_{2}}}\left(z-z_{l_{1}l_{2}}\right)^{\top}\hat{f}_{z}\left(z_{l_{1}l_{2}}\right)\hat{\ell}\left(z\right)dz+\frac{1}{2}\int_{R_{l_{1}l_{2}}}\left(z-z_{l_{1}l_{2}}\right)^{\top}\hat{f}_{zz}\left(\zeta_{l_{1}l_{2}}\left(z\right)\right)\left(z-z_{l_{1}l_{2}}\right)\hat{\ell}\left(z\right)dz\right).\label{eq: E^2 - formula 1}
\end{align}
For the first term in (\ref{eq: E^2 - formula 1}), we consider the
following expansion of $\hat{\ell}$ for $z\in R_{l_{1}l_{2}}$:
\[
\hat{\ell}\left(z\right)=\hat{\ell}\left(z_{l_{1}l_{2}}\right)+\left(z-z_{l_{1}l_{2}}\right)^{\top}\hat{\ell}_{z}\left(\eta_{l_{1}l_{2}}\left(z\right)\right),
\]
where $\eta_{l_{1}l_{2}}\left(z\right)\in R_{l_{1}l_{2}}$ depends
on $z$. Then,
\[
\int_{R_{l_{1}l_{2}}}\left(z-z_{l_{1}l_{2}}\right)^{\top}\hat{f}_{z}\left(z_{l_{1}l_{2}}\right)\hat{\ell}\left(z\right)dz=\sum_{i=1}^{2}\sum_{j=1}^{2}\int_{R_{l_{1}l_{2}}}\left(z-z_{l_{1}l_{2}}\right)^{\left(i\right)}\left(z-z_{l_{1}l_{2}}\right)^{\left(j\right)} \hat{f}_{z}^{\left(i\right)}\left(z_{l_{1}l_{2}}\right)\hat{\ell}_{z}^{\left(j\right)}\left(\eta_{l_{1}l_{2}}\left(z\right)\right)dz,
\]
where the term multiplying $\hat{\ell}\left(z_{l_{1}l_{2}}\right)$ vanishes due to symmetry. Consequently,
\begin{equation}
\left|\int_{R_{l_{1}l_{2}}}\left(z-z_{l_{1}l_{2}}\right)^{\top} \hat{f}_{z}\left(z_{l_{1}l_{2}}\right) \hat{\ell}\left(z\right)dz\right| \leq h_{z}^{2} \max_{1 \le i,j\le 2} \int_{R_{l_{1}l_{2}}} \left|\hat{f}_{z}^{\left(i\right)}\left(z_{l_{1}l_{2}}\right)\right| \left|\hat{\ell}_{z}^{\left(j\right)}\left(\eta_{l_{1}l_{2}}\left(z\right)\right)\right|dz.\label{eq: E^2 - formula 3}
\end{equation}
Next, for $j=1,2$ and any $w \in R_{l_{1}l_{2}}$, we have
\begin{align*}
\left|\hat{\ell}_{z}^{\left(j\right)}\left(w\right)\right| & =\left|2w^{\left(j\right)}\ell\left(w\right)+\left\Vert w\right\Vert ^{2}\ell_{z}^{\left(j\right)}\left(w\right)\right|\\
 & \leq2\left\Vert w\right\Vert \ell\left(w\right)+\left\Vert w\right\Vert ^{2}\left|\ell_{z}^{\left(j\right)}\left(w\right)\right|\\
 & \leq\left(2C_{\ell}\left(z_{\max}^{\mathbf{II}}\right)+C_{\ell}^{\prime}\left(z_{\max}^{\mathbf{II}}\right)\right)\left\Vert w\right\Vert ^{-1-A_{\ell}},
\end{align*}
by virtue of the upper bounds (\ref{eq:  intro - upper bound of the L=0000E9vy measure}) for the L{\'e}vy measure.
It is readily verified that $\left\Vert w\right\Vert \ge \frac{1}{3} \left\Vert z\right\Vert$ whenever $w, z\in R_{l_{1}l_{2}} \subset R_{z}^{\mathbf{II}}$.
Defining 
\[
\hat{C}_{\ell}=3^{1+A_{\ell}}\left(2C_{\ell}\left(z_{\max}^{\mathbf{II}}\right)+C_{\ell}^{\prime}\left(z_{\max}^{\mathbf{II}}\right)\right),
\]
it thus follows that
\begin{equation}
\left|\hat{\ell}_{z}^{\left(j\right)}\left(\eta_{l_{1}l_{2}}\left(z\right)\right)\right|\leq\hat{C}_{\ell}\left\Vert z\right\Vert ^{-1-A_{\ell}}
\label{eq: E^2 formula 4}
\end{equation}
for $j=1,2$. Combining (\ref{eq: proof - bounds on derivatives}), (\ref{eq: E^2 - formula 3}) and (\ref{eq: E^2 formula 4}) yields 
\begin{equation}
\left|\int_{R_{l_{1}l_{2}}}\left(z-z_{l_{1}l_{2}}\right)^{\top}\hat{f}_{z}\left(z_{l_{1}l_{2}}\right)\hat{\ell}\left(z\right)dz\right| \leq 3C_{f}^{\prime} \hat{C}_\ell h_{z}^{2}\left(\int_{R_{l_{1}l_{2}}}\left\Vert z\right\Vert ^{-2-A_{\ell}}dz\right).
\label{eq: E^2 - first-term bound}
\end{equation}
For the second term in (\ref{eq: E^2 - formula 1}) there holds
\begin{align*}
 & \left|\int_{R_{l_{1}l_{2}}}\left(z-z_{l_{1}l_{2}}\right)^{\top}\hat{f}_{zz}\left(\zeta_{l_{1}l_{2}}\left(z\right)\right)\left(z-z_{l_{1}l_{2}}\right)\hat{\ell}\left(z\right)dz\right|\nonumber \\
 & =\left|\sum_{i=1}^{2}\sum_{j=1}^{2}\int_{R_{l_{1}l_{2}}}\left(z-z_{l_{1}l_{2}}\right)^{\left(i\right)}\left(z-z_{l_{1}l_{2}}\right)^{\left(j\right)} \hat{f}_{zz}^{\left(i,j\right)}\left(\zeta_{l_{1}l_{2}}\left(z\right)\right) \hat{\ell}\left(z\right)dz\right|\nonumber \\
 & \leq h_{z}^{2} \max_{1 \le i,j\le 2} \int_{R_{l_{1}l_{2}}} \left|\hat{f}_{zz}^{\left(i,j\right)}\left(\zeta_{l_{1}l_{2}}\left(z\right)\right)\right| \hat{\ell}\left(z\right)dz.
\end{align*}
Invoking (\ref{eq: proof - bounds on derivatives}) and (\ref{eq:  intro - upper bound of the L=0000E9vy measure}), it follows that
\begin{equation}
\left|\int_{R_{l_{1}l_{2}}}\left(z-z_{l_{1}l_{2}}\right)^{\top}\hat{f}_{zz}\left(\zeta_{l_{1}l_{2}}\left(z\right)\right)\left(z-z_{l_{1}l_{2}}\right)\hat{\ell}\left(z\right)dz\right| \le
3C_{f}^{\prime\prime} C_{\ell}\left(z_{\max}^{\mathbf{II}}\right) h_{z}^{2} \left(\int_{R_{l_{1}l_{2}}}\left\Vert z\right\Vert ^{-2-A_{\ell}}dz\right).
\label{eq: E^2 - second-term bound}
\end{equation}
By (\ref{eq: E^2 - first-term bound}) and  (\ref{eq: E^2 - second-term bound}),
the error $E^{\left(2\right)}$ satisfies the following bound 
\[
\left|E^{\left(2\right)}\right|\leq C\left(\int_{R_{z}^{\mathbf{II}}}\left\Vert z\right\Vert ^{-2-A_{\ell}}dz\right)h_{z}^{2}
\]
with constant
\[
C = 3C_{f}^{\prime} \hat{C}_\ell + \frac{3}{2}C_{f}^{\prime\prime} C_{\ell}\left(z_{\max}^{\mathbf{II}}\right).
\]
Using (\ref{eq: proof - radial integral}) and (\ref{eq: proof - radial integral added}), we obtain that 
\[
\int_{R_{z}^{\mathbf{II}}}\left\Vert z\right\Vert ^{-2-A_{\ell}}dz
=\begin{cases}
O \Big( \left(z_{\max}^{\mathbf{II}}\right)^{-A_{\ell}}+h_{z}^{-A_{\ell}}\Big) & \textrm{if } 0<A_{\ell}<2,\\
O \Big(\ln (\sqrt{2}\,z_{\max}^{\mathbf{II}})-\ln (Mh_{z}) \Big) & \textrm{if } A_{\ell}=0.
\end{cases}
\]
Since for any given $\epsilon>0$ it holds that $\ln h = O\left(h^{-\epsilon}\right)$ as $h\rightarrow0^{+}$, there follows 
\[
\int_{R_{z}^{\mathbf{II}}}\left\Vert z\right\Vert ^{-2-A_{\ell}}dz
=\begin{cases}
O\left(h_{z}^{-A_{\ell}}\right) & \textrm{if } 0<A_{\ell}<2,\\
O\left(h_{z}^{-\epsilon}\right) & \textrm{if } A_{\ell}=0.
\end{cases}
\]
Thus, 
\begin{equation}
E^{\left(2\right)} = \begin{cases}
O\left(h_{z}^{2-A_{\ell}}\right) & \textrm{if } 0<A_{\ell}<2,\\
O\left(h_{z}^{2-\epsilon}\right) & \textrm{if } A_{\ell}=0,
\end{cases}
\label{eq: E^2}
\end{equation}
uniformly in $\left(x,t\right)\in S$.

We conclude with $E^{\left(3\right)}$. It is well known, see e.g.
\citet{quarteroni2007}, that the composite midpoint rule has second-order convergence for smooth integrands. Since $z_{\max}^{\mathbf{II}}$ and $z_{\max}^{\mathbf{III}}$ are independent of $h_{z}$ and $f\cdot\ell\in\mathcal{C}^{2}\left(R_{z}^{\mathbf{III}}\right)$,
there holds 
\begin{equation}
E^{\left(3\right)}=O\left(h_{z}^{2}\right)\label{eq: E^3}
\end{equation}
uniformly in $\left(x,t\right)\in S$.

From (\ref{eq: errore decomposition}), (\ref{eq: E^1}), (\ref{eq: E^2}) and (\ref{eq: E^3}), the stated result follows.

\section{\label{sec: appendix - summation operator as a circulant matrix-vector multiplication}Summation
operator as a circulant matrix-vector multiplication}

Let $N_{z}=1$, $N_{y}^{-}=0$ and $N_{y}^{+}=1$. Then $\sharp\text{out}=N_{y}^{+}+N_{y}^{-}+1=2$ and $\sharp\text{in}=2N_{z}+N_{y}^{+}+N_{y}^{-}=3$. The quadrature matrix $\Omega$, whose entries are the coefficients $\omega$ defined in \eqref{eq: integral disc - omega coefficients}, is given by 
\[
\Omega=\left[\begin{array}{cc}
\omega_{-1,-1} & \omega_{-1,0}\\
\omega_{0,-1} & \omega_{0,0}
\end{array}\right]\in\mathbb{R}^{2N_{z}\times2N_{z}}.
\]
The first row of the circulant matrix $C$ is defined according to 
\[
C_{1,\cdot}=\text{vec}\left(\left[\begin{array}{ccc}
\omega_{-1,-1} & \omega_{-1,0} & 0\\
\omega_{0,-1} & \omega_{0,0} & 0\\
0 & 0 & 0
\end{array}\right]\right)\in\mathbb{R}^{\left(\sharp\text{in}\right)^{2}\times1},
\]
while the entire matrix is
\[
C=\left[\begin{array}{ccccccccc}
{\color{red}\omega_{-1,-1}} & {\color{red}\omega_{0,-1}} & {\color{red}0} & {\color{red}\omega_{-1,0}} & {\color{red}\omega_{0,0}} & {\color{red}0} & {\color{red}0} & {\color{red}0} & {\color{red}0}\\
{\color{red}0} & {\color{red}\omega_{-1,-1}} & {\color{red}\omega_{0,-1}} & {\color{red}0} & {\color{red}\omega_{-1,0}} & {\color{red}\omega_{0,0}} & {\color{red}0} & {\color{red}0} & {\color{red}0}\\
0 & 0 & \omega_{-1,-1} & \omega_{0,-1} & 0 & \omega_{-1,0} & \omega_{0,0} & 0 & 0\\
{\color{red}0} & {\color{red}0} & {\color{red}0} & {\color{red}\omega_{-1,-1}} & {\color{red}\omega_{0,-1}} & {\color{red}0} & {\color{red}\omega_{-1,0}} & {\color{red}\omega_{0,0}} & {\color{red}0}\\
{\color{red}0} & {\color{red}0} & {\color{red}0} & {\color{red}0} & {\color{red}\omega_{-1,-1}} & {\color{red}\omega_{0,-1}} & {\color{red}0} & {\color{red}\omega_{-1,0}} & {\color{red}\omega_{0,0}}\\
\omega_{0,0} & 0 & 0 & 0 & 0 & \omega_{-1,-1} & \omega_{0,-1} & 0 & \omega_{-1,0}\\
\omega_{-1,0} & \omega_{0,0} & 0 & 0 & 0 & 0 & \omega_{-1,-1} & \omega_{0,-1} & 0\\
0 & \omega_{-1,0} & \omega_{0,0} & 0 & 0 & 0 & 0 & \omega_{-1,-1} & \omega_{0,-1}\\
\omega_{0,-1} & 0 & \omega_{-1,0} & \omega_{0,0} & 0 & 0 & 0 & 0 & \omega_{-1,-1}
\end{array}\right]\in\mathbb{R}^{\left(\sharp\text{in}\right)^{2}\times\left(\sharp\text{in}\right)^{2}}.
\]
The entries in the first, second, fourth, and fifth rows (highlighted in red) correspond to the matrix $\tilde{I}C\in\mathbb{R}^{\left(\sharp\text{out}\right)^{2}\times\left(\sharp\text{in}\right)^{2}}$.

\newpage{}

\bibliographystyle{zootaxa}
\phantomsection\addcontentsline{toc}{section}{\refname}\bibliography{bibliography}

\end{document}